\documentclass[a4paper, 12pt, draft]{article}
\makeatletter
\@setfontsize\normalsize\@xiipt{15.5}%
\renewcommand*{\author}[1]{\gdef\@author{#1}\gdef\@pauthor{{\def\and{ --- }#1}}}
\renewcommand*{\title}[1]{\gdef\@title{#1}\gdef\@ptitle{#1}}
\if@twoside         
  \def\ps@draft{%
    \def\@oddfoot{\small\null\hfil\thepage\hfil}
    \let\@evenfoot\@oddfoot
    \def\@evenhead{\small\@date\hfil\slshape\@pauthor\hfil}
    \def\@oddhead{\small\null\hfil\slshape\@ptitle\hfil}
    \let\@mkboth\@gobbletwo
    \let\sectionmark\@gobble
    \let\subsectionmark\@gobble
   }
\else               
  \def\ps@draft{%
    \def\@oddfoot{\small\@date\hfil\slshape\@pauthor\hfil\upshape\thepage}
    \def\@oddhead{\small\null\hfil\slshape\@ptitle\hfil}
    \let\@mkboth\@gobbletwo
    \let\sectionmark\@gobble
    \let\subsectionmark\@gobble
  }
\fi

\newcommand*{\keywords}[1]{\gdef\@keywords{#1}}
\newcommand{\keywordsname}{Key words and phrases}
\newcommand*{\subjclass}[1]{\gdef\@subjclass{#1}}
\newcommand{\subjclassname}{1991 AMS Mathematics Subject Classification}
\def\@maketitle{%
  \newpage
  \null
  \vskip 2em%
  \begin{center}%
    {\Large\bfseries \@title \par}%
    \vskip 1.5em%
    {\small\scshape
      \lineskip .5em%
      \begin{tabular}[t]{c}%
        \@author
      \end{tabular}\par
    }%
    \vskip 1em%
    {\small\@date}
  \end{center}%
  \par
  \vskip 1.5em
  \begingroup
    \let\@makefnmark\relax \let\@thefnmark\relax
    \ifx\@empty\@subjclass\else
       \@footnotetext{{\itshape\subjclassname}.\enspace\@subjclass.}
    \fi
    \ifx\@empty\@keywords\else
       \@footnotetext{{\itshape\keywordsname}.\enspace\@keywords.}
    \fi
  \endgroup
}
\renewcommand{\part}{\par
   \addvspace{4ex}%
   \@afterindentfalse
   \secdef\@part\@spart}

\def\@part[#1]#2{%
    \ifnum \c@secnumdepth >\m@ne
      \refstepcounter{part}%
      \addcontentsline{toc}{part}{\thepart\hspace{1em}#1}%
    \else
      \addcontentsline{toc}{part}{#1}%
    \fi
    {\parindent \z@ \raggedright
     \interlinepenalty \@M
     \reset@font
     \ifnum \c@secnumdepth >\m@ne
       \large\bfseries \partname~\thepart
       \par\nobreak
     \fi
     \Large \bfseries #2%
     \markboth{}{}\par}%
    \nobreak
    \vskip 3ex
    \@afterheading}
\def\@spart#1{%
    {\parindent \z@ \raggedright
     \interlinepenalty \@M
     \reset@font
     \Large \bfseries #1\par}%
     \nobreak
     \vskip 3ex
     \@afterheading}
\def\@endpart{\vfil\newpage
              \if@twoside
                \hbox{}%
                \thispagestyle{empty}%
                \newpage
              \fi
              \if@tempswa
                \twocolumn
              \fi}
\renewcommand{\section}{\@startsection {section}{1}{\z@}%
                                   {-3.5ex \@plus -1ex \@minus -.2ex}%
                                   {2.3ex \@plus.2ex}%
                                   {\reset@font\large\bfseries}}
\renewcommand{\subsection}{\@startsection{subsection}{2}{\z@}%
                                     {-3.25ex\@plus -1ex \@minus -.2ex}%
                                     {1.5ex \@plus .2ex}%
                                     {\reset@font\normalsize\bfseries}}
\renewcommand{\subsubsection}{\@startsection{subsubsection}{3}{\z@}%
                                     {-3.25ex\@plus -1ex \@minus -.2ex}%
                                     {1.5ex \@plus .2ex}%
                                     {\reset@font\normalsize\bfseries}}
\renewcommand{\paragraph}{\@startsection{paragraph}{4}{\z@}%
                                    {3.25ex \@plus1ex \@minus.2ex}%
                                    {-1em}%
                                    {\reset@font\normalsize\bfseries}}
\renewcommand{\subparagraph}{\@startsection{subparagraph}{5}{\parindent}%
                                       {3.25ex \@plus1ex \@minus .2ex}%
                                       {-1em}%
                                      {\reset@font\normalsize\bfseries}}

\renewcommand{\theenumi}{\alph{enumi}}

\renewcommand{\theenumii}{\roman{enumii}}

\renewcommand{\p@enumii}{\theenumi.}
\renewcommand{\theenumiii}{\Alph{enumiii}}

\renewcommand{\p@enumiii}{\theenumi.\theenumii.}

\renewcommand{\p@enumiv}{\p@enumiii\theenumiii.}
\@addtoreset{equation}{section}

\RequirePackage{amsthm}

\renewenvironment{proof}[1][\proofname]{\par
  \normalfont
  \topsep6\p@\@plus6\p@ \trivlist
  \item[\hskip\labelsep\slshape
    #1\@addpunct{.}]\ignorespaces
}{%
  \qed\endtrivlist
}

\theoremstyle{plain}

\newtheorem{theorem}{Theorem}[section]
\newtheorem{proposition}[theorem]{Proposition}
\newtheorem{lemma}[theorem]{Lemma}
\newtheorem{corollary}[theorem]{Corollary}

\theoremstyle{definition}

\newtheorem{definition}[theorem]{Definition}

\newtheorem{remark}[theorem]{Remark}

\RequirePackage{amsmath}
\RequirePackage{amssymb}

\makeatother

\font\goth=eufm10 at 12pt
\font\subgoth=eufm10 at 10pt

\newcommand{\gothg}{\mbox{\goth g}}

\newcommand{\gothh}{\mbox{\goth h}}
\newcommand{\gothp}{\mbox{\goth p}}

\newcommand{\gothl}{\mbox{\goth l}}

\newcommand{\gothn}{\mbox{\goth n}}
\newcommand{\gothz}{\mbox{\goth z}}

\newcommand{\sgothp}{\mbox{\subgoth p}}

\newcommand{\sgothl}{\mbox{\subgoth l}}

\def\shb{\mathcal{B}}
\def\shc{\mathcal{C}}

\def\shf{\mathcal{F}}

\def\shl{\mathcal{L}}
\def\shm{\mathcal{M}}
\def\shn{\mathcal{N}}

\def\shu{\mathcal{U}}

\newcommand{\C}{\mathbb{C}}

\newcommand{\Z}{\mathbb{Z}}
\DeclareMathOperator{\codim}{codim}

\DeclareMathOperator{\id}{id}
\newcommand{\BDC}{\mathbf{D}^{\mathrm{b}}}

\renewcommand{\to}{\rightarrow}
\newcommand{\from}{\leftarrow}
\newcommand{\Ker}{{\rm Ker}}
\newcommand{\isoto}{\rightarrow^{\!\!\!\!\!\!\!\sim}}

\newcommand{\Hom}[1][]{\mathrm{Hom}_{\raise1.5ex\hbox to.1em{}#1}}
\newcommand{\RHom}[1][]{\mathrm{RHom}_{\raise1.5ex\hbox to.1em{}#1}}
\newcommand{\THom}[1][]{\mathrm{THom}_{\raise1.5ex\hbox to.1em{}#1}}
\newcommand{\Mod}{\mathrm{Mod}}
\newcommand{\Tens}[1][]{\mathbin{\otimes_{\raise1.5ex\hbox to-.1em{}#1}}}
\newcommand{\LTens}[1][]{\mathbin{\otimes_{\raise1.5ex\hbox to-.1em{}#1}^{L}}}
\newcommand{\sect}{\Gamma}
\newcommand{\rsect}{\mathrm{R}\Gamma}
\renewcommand{\hom}[1][]{{\mathcal{H}om}_{\raise1.5ex\hbox to.1em{}#1}}
\newcommand{\rhom}[1][]{{R\mathcal{H}om}_{\raise1.5ex\hbox to.1em{}#1}}
\newcommand{\thom}[1][]{{T\mathcal{H}om}_{\raise1.5ex\hbox to.1em{}#1}}
\newcommand{\tens}[1][]{\mathbin{\otimes_{\raise1.5ex\hbox to-.1em{}#1}}}
\newcommand{\ltens}[1][]{\mathbin{\otimes_{\raise1.5ex\hbox to-.1em{}#1}^{L}}}

\newcommand\etens{\mathbin{\boxtimes}}

\newcommand{\roim}[1]{{R#1}_*}
\newcommand{\reim}[1]{{R#1}_!}
\newcommand{\opb}[1]{#1^{-1}}

\newcommand{\D}{\mathcal{D}}
\renewcommand{\O}{\mathcal{O}}

\newcommand{\B}{\mathcal{B}}
\newcommand{\I}{\mathcal{I}}

\newcommand{\doim}[1]{\underline{#1}_*}
\newcommand{\deim}[1]{\underline{#1}_!}
\newcommand{\dopb}[1]{\underline{#1}^{-1}}
\renewcommand{\deim}[1]{\underline{f}_!}
\newcommand{\End}[1]{\mbox{\rm End}_{#1}}
\newcommand{\eps}{\varepsilon}
 
\newcommand{\Ad}{{\rm Ad}} 
\newcommand{\lam}{{\lambda}} 
\newcommand{\an}{{\rm an}}

\newcommand{\maprightu}[1]{%
\smash{\mathop{%
\hbox to 1cm{\rightarrowfill}}\limits^{#1}}}
\newcommand{\maprightd}[1]{%
\smash{\mathop{%
\hbox to 1cm{\rightarrowfill}}\limits_{#1}}}
\newcommand{\mapleftu}[1]{%
\smash{\mathop{%
\hbox to 1cm{\leftarrowfill}}\limits^{#1}}}
\newcommand{\mapleftd}[1]{%
\smash{\mathop{%
\hbox to 1cm{\leftarrowfill}}\limits_{#1}}}

\newcommand{\downeq}{\vline\kern2pt\hbox{$\wr$}}

\title{Radon transforms for quasi-equivariant $\D$-modules on generalized flag manifolds}
\author{Corrado Marastoni \and Toshiyuki Tanisaki}
\begin{document}

\maketitle

\begin{abstract}
In this paper we deal with Radon transforms for generalized flag manifolds in the framework of quasi-equivariant $\D$-modules.
We shall follow the method employed by Baston-Eastwood and analyze the Radon transform using the Bernstein-Gelfand-Gelfand resolution and the Borel-Weil-Bott theorem.
We shall determine the transform completely on the level of the Grothendieck groups.
Moreover, we point out a vanishing criterion and give a sufficient condition in order that a $\D$-module associated to an equivariant locally free $\O$-module is transformed into  an object of the same type.
The case of maximal parabolic subgroups of classical simple groups is studied in detail. 
\end{abstract}

\section*{Introduction} 
\addcontentsline{toc}{section}{Introduction}

Let $G$ be a reductive algebraic group over $\C$, $P$ and $Q$ two parabolic
subgroups containing the same Borel subgroup of $G$. Let $X=G/P$, $Y=G/Q$, and let $S$ be the unique closed $G$-orbit in $X\times Y$ for the diagonal action.
Then we can identify $S$ with $G/P\cap Q$.
The natural correspondence
$$X\overset{f}{\longleftarrow} S\overset{g}{\longrightarrow} Y,$$
where $f$ and $g$ are the restriction to $S$ of the projections of $X\times Y$
on $X$ and $Y$, induces an integral transform from $X$ to $Y$ which generalizes
the classical Radon-Penrose transform. 
This subject has been investigated intensively both
in the complex and real domains (see e.g.\ 
Baston-Eastwood \cite{B-E},
D'Agnolo-Schapira \cite{D'A-S}, 
Kakehi \cite{Kakehi},
Marastoni \cite{Ma}, 
Oshima \cite{Osh}, 
Sekiguchi \cite{Sekiguchi}, 
Tanisaki \cite{Ta}). 

Our aim is to study this transform in the framework of quasi-$G$-equivariant
$\D$-modules (see Kashiwara~\cite{Kas}), i.e.\ the functor
\begin{equation}
\label{eq:RadonDModIntro}
R:\BDC_G(\D_X)\to \BDC_G(\D_Y), \qquad R(\shm)=\doim{g}\dopb{f}\shm,
\end{equation}
where $\BDC_G(\D_.)$ denotes the derived category of quasi-G-equivariant
$\D$-modules with bounded cohomologies, and $\doim{g}$ and $\dopb{f}$ are the 
operations of direct image (integration) and inverse image (pull-back) for
$\D$-modules. 
More precisely, we consider a $\D_X$-module of type
$\shm=\D\shl=\D_X\tens_{\O_X}\shl$, where $\shl$ is an irreducible $G$-equivariant locally free $\O_X$-module. 
In this case it is easily seen that
\begin{equation}
\label{eq:intro0}
H^p(R(\D\shl))=0
\quad
\mbox{for any $p<0$}
\end{equation}
(see Lemma~\ref{lem:lower vanishing} below).
Note  that the Grothendieck group of the category of quasi-$G$-equivariant $\D_X$-modules of finite length is spanned by elements corresponding to the objects of the form $\D\shl$.

As in Baston-Eastwood~\cite{B-E} our analysis relies on the Bernstein-Gelfand-Gelfand resolution and the Borel-Weil-Bott theorem.
Using the Bernstein-Gelfand-Gelfand resolution in the parabolic setting (see Bernstein-Gelfand-Gelfand~\cite{BGG}, Lepowsky~\cite{Lep}, Rocha-Caridi~\cite{Rocha-Caridi}) we obtain a resolution of the quasi-$G$-equivariant $\D_S$-module $\dopb{f}(\D\shl)$ of the form:
\begin{equation}
\label{eq:intro1}
0\to\bigoplus_{k=1}^{r_n}\D\shl_{nk}
\to\cdots
\to\bigoplus_{k=1}^{r_0}\D\shl_{0k}
\to\dopb{f}(\D\shl)\to0,
\end{equation}
where $\shl_{ik}$ are irreducible $G$-equivariant locally free $\O_S$-modules (see \S~\ref{subsection:RTQED} for the explicit description of $\shl_{ik}$).
Then we have 
$$
\doim{g}(\D\shl_{ik})=\D_Y\tens[\O_Y]Rg_*(\shl_{ik}\tens[\O_S]\Omega_g)
$$ 
by the definition of $\doim{g}$, where $\Omega_g$ denotes the sheaf of relative differential forms with maximal degree along the fibers of $g$.
Moreover, the Borel-Weil-Bott theorem tells us the structure of $Rg_*(\shl_{ik}\tens[\O_S]\Omega_g)$.
In particular, we have either $Rg_*(\shl_{ik}\tens[\O_S]\Omega_g)=0$ or there exist a non-negative integer $m_{ik}$ and an irreducible $G$-equivariant $\O_Y$-module $\shl'_{ik}$ such that $Rg_*(\shl_{ik}\tens[\O_S]\Omega_g)=\shl'_{ik}[-m_{ik}]$.
Thus setting
$$
\I=\{(i,k)\,;\,0\leqq i\leqq n, 1\leqq k\leqq r_i, Rg_*(\shl_{ik}\tens[\O_S]\Omega_g)\ne0\},
$$
we have
\begin{equation}
\label{eq:intro2}
\doim{g}(\D\shl_{ik})=
\left\{
\begin{array}{ll}
\D\shl'_{ik}[-m_{ik}]\qquad
&{((i,k)\in\I)},\\
0
&{((i,k)\not\in\I)}
\end{array}
\right.
\end{equation}
(see \S\ref{subsection:RTQED} below for concrete descriptions of $\I$ and $\shl_{ik}, m_{ik}$ for $(i,k)\in \I$).

Then we can study the structure of $R(\D\shl)=\doim{g}\dopb{f}(\D\shl)$ using \eqref{eq:intro0}, \eqref{eq:intro1} and \eqref{eq:intro2}.
For example we have the following result.
\begin{theorem}
\label{thm:intro1}
Let the notation be as above.
\begin{itemize}
\item[\rm(i)]
We have
$$
\sum_p(-1)^p[H^p(R(\D\shl))]=\sum_{(i,k)\in\I}(-1)^{i-m_{ik}}[\D\shl'_{ik}]
$$
in the Grothendieck group of the category of quasi-$G$-equivariant $\D_Y$-modules.
\item[\rm(ii)]
If $\I=\emptyset$, then $R(\D\shl)=0$.
\item[\rm(iii)]
If $\I$ consists of a single element $(i,k)$, then $R(\D\shl)=\D\shl'_{ik}[i-m_{ik}]$.
\item[\rm(iv)]
If $i\geqq m_{ik}$ for any $(i,k)\in\I$, then we have $H^p(R(\D\shl))=0$ unless $p=0$.
\item[\rm(v)]
If $i> m_{ik}$ for any $(i,k)\in\I$ with $i>0$ and if $m_{01}=0$, then there exists an epimorphism $\D_Y\shl'_{01}\to H^0(R(\D\shl))$ $($note that $r_0=1$ $)$.
\end{itemize}
\end{theorem}

Assume that $\shl$ is invertible and that there exists a $G$-equivariant invertible  $\O_Y$-module $\shl'$ satisfying $f^*\shl\tens[\O_S]\Omega_g=g^*\shl'$.
We call such a pair $(\shl,\shl')$ an extremal case for the correspondence (if $P\cup Q$ generates the group $G$ and if $G$ is semisimple, then there exists a unique extremal case).
In this case there exists a natural nontrivial $\D_Y$-linear morphism
\begin{equation}
\Phi:\D\shl'\ \to H^0(R(\D\shl)).
\end{equation}
\begin{theorem}
\label{thm:intro2}
Let $(\shl,\shl')$ be an extremal case.
\begin{itemize}
\item[\rm(i)]
We have $H^p(R(\D\shl))=0$ for any $p\ne0$ if and only if $i\geqq m_{ik}$ for any $(i,k)\in\I$.
\item[\rm(ii)]
Assume that $H^p(R(\D\shl))=0$ for any $p\ne0$.
Then $\Phi$ is an epimorphism if and only if $i> m_{ik}$ for any $(i,k)\in\I$ with $i>0$.
\item[\rm(iii)]
Assume that $H^p(R(\D\shl))=0$ for any $p\ne0$.
Then $\Phi$ is an isomorphism if and only if $\I$ consists of a single element $(0,1)$.
\end{itemize}
\end{theorem}
We do not know an example of an extremal case $(\shl,\shl')$ such that $H^p(R(\D\shl))\ne0$ for some $p\ne0$.
We have checked that $H^p(R(\D\shl))=0$ for any $p\ne0$ by a case-by-case analysis when $G$ is a classical simple group, $P$, $Q$ are maximal parabolic subgroups and $(\shl,\shl')$ is the extremal case.
In general the morphism $\Phi$ for an extremal case $(\shl,\shl')$ is not necessarily an epimorphism nor a monomorphism.
It would be an interesting problem to determine the kernel and the cokernel of $\Phi$.

The transform of a $\D$-module, a problem of analytic nature, is not sufficient to cover the general problem of integral geometry.
In order to do this, one should couple the transforms in the frameworks of $\D$-modules and sheaves. 
This is better described in the adjunction formulas (see D'Agnolo-Schapira \cite{D'A-S}), and we shall briefly discuss this point with an example in the case of $G=SL_{n+1}(\C)$.

\vskip .3cm
We would like to thank M.\ Kashiwara for useful conversation on quasi-equivariant $\D$-modules.

\section{Preliminaries on $\D$-modules}
\subsection{Functors for $\D$-modules}

Let $Z$ be an algebraic manifold (smooth algebraic variety) over $\C$. 
We denote by $\O_Z$ the structure sheaf, by $\Omega_Z$ the invertible $\O_Z$-module of differential forms of maximal degree, and by $\D_Z$ the sheaf of differential operators.
In this paper an $\O_Z$-module means a quasi-coherent $\O_Z$-module and a $\D_Z$-module means a left $\D_Z$-module which is quasi-coherent over $\O_Z$.
We denote by $\Mod(\D_Z)$ the category of $\D_Z$-modules and by $\BDC(\D_Z)$ the derived category of $\Mod(\D_Z)$ whose objects have bounded cohomology. 

If $f:Z\to Z'$ is a morphism, we set
$$
\Omega_f=\Omega_{Z/Z'}=\Omega_Z\tens[\opb{f}\O_{Z'}]\opb{f}\Omega_{Z'}^{\tens-1};
$$
and, for an $\O_{Z'}$-module $\shl'$, we set
$$
f^*\shl'=\O_Z\tens[\opb{f}\O_{Z'}]\opb{f}\shl',\qquad
Lf^*\shl'=\O_Z\ltens[\opb{f}\O_{Z'}]\opb{f}\shl'.
$$
We denote by $\doim{f}$ and $\dopb{f}$ the direct and inverse image
for left $\D$-modules:
{\small
$$
\begin{array}{ll}
\doim{f}:\BDC(\D_Z)\to \BDC(\D_{Z'}), & \doim{f}\shm=\roim{f}(\D_{Z'\from
Z}\ltens[\D_Z]\shm), \\
\dopb{f}:\BDC(\D_{Z'})\to \BDC(\D_{Z}), &
\dopb{f}\shm'=\D_{Z\to Z'}\ltens[\opb{f}\D_{Z'}]\opb{f}\shm',
\end{array}
$$
}where a $(\D_Z,f^{-1}\D_{Z'})$-bimodule $\D_{Z\to Z'}$ and an $(f^{-1}\D_{Z'},\D_Z)$-bimodule $\D_{Z'\from Z}$ are defined by
$$
\D_{Z\to Z'}=\O_{Z}\tens[\opb{f}\O_{Z'}]\opb{f}\D_{Z'},\qquad
\D_{Z'\from Z}=
\Omega_Z\tens[\O_Z]\D_{Z\to Z'}\tens[\opb{f}\O_{Z'}]\opb{f}\Omega_{Z'}^{\tens-1}.
$$
Note that for a $\D_{Z'}$-module $\shm$ we have $\dopb{f}\shm\simeq Lf^*\shm$ as a complex of $\O_Z$-modules.
Note also that we have canonical morphisms $\O_Z\to\D_{Z\to Z'}$ and $\Omega_f\to\D_{Z'\from Z}$ of $\O_Z$-modules.

The following result is well-known and easy to prove.
\begin{lemma}
\label{lem:transfer}
Let $f_1:Z\to X_1$ and $f_2:Z\to X_2$ be morphisms of algebraic manifolds.
\begin{itemize}
\item[\rm(i)]
We have
$$
\D_{X_2\from Z}\ltens[\D_Z]\D_{Z\to X_1}\isoto
\opb{f_1}\Omega_{X_1}\ltens[\opb{f_1}\O_{X_1}]
(\D_{X_1\times X_2\from Z}\ltens[\D_Z]\O_Z).
$$
\item[\rm(ii)]
Assume that $Z\to X_1\times X_2$ is an embedding. 
Then we have 
$$\D_{X_2\from Z}\ltens[\D_Z]\D_{Z\to X_1}
=\D_{X_2\from Z}\tens[\D_Z]\D_{Z\to X_1},$$
and the canonical morphism $\Omega_{f_2}\to\D_{X_2\from Z}\tens[\D_Z]\D_{Z\to X_1}$ of $(\opb{f_2}\O_{X_2},\opb{f_1}\O_{X_1})$-bimodules is a monomorphism.
\end{itemize}
\end{lemma}
For a locally free $\O_{Z}$-module
$\shl$, we set
$$
\D\shl=\D_{Z}\tens[\O_{Z}]\shl,
$$
and for a closed submanifold $Z$ of an algebraic manifold $X$ we define a $\D_X$-module $\B_{Z|X}$ supported on $Z$ by
$$
\B_{Z|X}=H^{d}_{[Z]}(\O_{X})
=\doim{i}\O_Z,
$$
where $d=\codim_{X}Z$ and $i:Z\to X$ denotes the embedding.
\subsection{Radon transforms}
Let $X$ and $Y$ be algebraic manifolds over $\C$, and denote by $q_1$ and 
$q_2$ the projections of $X\times Y$ onto $X$ and $Y$ respectively. 
Let $S$ be a locally closed submanifold of $X\times Y$ and let $i:S\to X\times Y$ be the embedding. 
The geometric correspondence
\begin{equation}
\label{eq:corr}
X\overset{f}{\longleftarrow} S\overset{g}{\longrightarrow}Y
\end{equation}
where $f$ and $g$ are the restrictions of $q_1$ and $q_2$, induces a functor
{\small
\begin{equation}
R:\BDC(\D_X)\to \BDC(\D_Y), \label{eq:Radon} 
\qquad
R(\shm)=\doim{g}\dopb{f}(\shm),
\end{equation}
}
called the Radon transform.
\begin{lemma}
\label{lem:RM}
Let $\shm$ be a $\D_X$-module.
\begin{itemize}
\item[\rm(i)]
We have
\begin{align*}
R(\shm)&=
Rg_*((\D_{Y\from S}\tens[\D_S]\D_{S\to X})\ltens[f^{-1}\D_X]f^{-1}\shm)\\
&=Rg_*(f^{-1}(\Omega_X\tens[\O_X]\shm)\ltens[f^{-1}\D_X]
(\D_{X\times Y\from S}\tens[\D_S]\O_S)).
\end{align*}
\item[\rm(ii)]
If $S$ is closed in $X\times Y$, then we have
$$
R(\shm)=\doim{q_2}(\dopb{q_1}\shm\ltens[\O_{X\times Y}]\B_{S|X\times Y}).
$$
\end{itemize}
\end{lemma}
\begin{proof}
(i) follows from the definition and Lemma~\ref{lem:transfer}, and (ii) is a consequence of the projection formula for $\D$-modules.
\end{proof}

Let us consider the special case where
$\shm=\D\shl=\D_X\tens[\O_X]\shl$.
By Lemma~\ref{lem:RM} we have the following.
\begin{lemma}
\label{lem:RDL}
Let $\shl$ be a locally free $\O_X$-module.
\begin{itemize}
\item[\rm(i)]
We have
\begin{align*}
R(\D\shl)&=
Rg_*((\D_{Y\from S}\tens[\D_S]\D_{S\to X})\tens[f^{-1}\O_X]f^{-1}\shl)\\
&=Rg_*(f^{-1}(\Omega_X\tens[\O_X]\shl)\tens[f^{-1}\O_X]
(\D_{X\times Y\from S}\tens[\D_S]\O_S)).
\end{align*}
\item[\rm(ii)]
If $S$ is closed in $X\times Y$, then we have
$$
R(\D\shl)=
Rq_{2*}(\opb{q_1}(\Omega_X\tens[\O_X]\shl)\tens[\opb{q_1}\O_X]\B_{S|X\times Y}).
$$
\end{itemize}
\end{lemma}
An immediate consequence of Lemma \ref{lem:RDL}(i) is:
\begin{lemma}
\label{lem:lower vanishing}
For any locally free $\O_X$-module $\shl$ we have
$H^p(R(\D\shl))=0$ for any $p<0$ .
\end{lemma}
\begin{definition}
\label{de:extrcase}
Let $\shl$ (resp.\ $\shl'$) be a locally free $\O_X$- (resp.\ $\O_Y$-)module.
We say that the pair $(\shl,\shl')$ is an extremal case for the
correspondence \eqref{eq:corr} if there is an $\O_S$-linear isomorphism
$$
\Omega_g\tens[\opb{f}\O_X]\opb{f}\shl \simeq {g}^*\shl'.
$$
\end{definition}
\begin{proposition}
\label{pr:ecnontrmor}
Let $(\shl,\shl')$ be an extremal case for \eqref{eq:corr}.
Then there exists a natural nontrivial $\D_Y$-linear morphism
\begin{equation}
\label{eq:Dmor}
\D\shl'\to H^0(R(\D\shl)).
\end{equation}
\end{proposition}
\begin{proof}
The canonical morphism
$\Omega_{g}\to\D_{Y\from S}\tens[\D_S]\D_{S\to X}$
 induces a monomorphism
$$
g^*\shl'\simeq\Omega_g\tens[\opb{f}\O_X]\opb{f}\shl\to
\D_{Y\from S}\tens[\D_S]\D_{S\to X}
\tens[\opb{f}\O_X]\opb{f}\shl
$$
of $\opb{g}\O_Y$-modules.
Applying $g_*$ we obtain a sequence of morphisms
\begin{eqnarray*}
\shl'&\to&\shl'\tens[\O_Y]g_*\O_S
\simeq
g_*(g^*\shl')\\
&\to&
g_*(\D_{Y\from S}\tens[\D_S]\D_{S\to X}
\tens[\opb{f}\O_X]\opb{f}\shl)
=
H^0(R(\D\shl))
\end{eqnarray*}
of $\O_Y$-modules.
The morphism $\shl'\to\shl'\tens[\O_Y]g_*\O_S$ is nontrivial by the definition, and the morphism
$g_*(g^*\shl')\to g_*(\D_{Y\from S}\tens[\D_S]\D_{S\to X}\tens[\opb{f}\O_X]\opb{f}\shl)$
is a monomorphism by the left exactness of $g_*$.
Thus the composition $\shl'\to H^0(R(\D\shl))$ is nontrivial.
Hence it induces a canonical nontrivial morphism $\D\shl'\to H^0(R(\D\shl))$ of $\D_Y$-modules.
\end{proof}

\subsection{Adjunction formulas}
In this subsection we consider topological problems, and hence we work in the analytic category rather than the algebraic category.

For a complex manifold $Z$ we denote by $\O_Z$ the sheaf of holomorphic functions on $Z$ and by $\D_Z$ the sheaf of holomorphic differential operators.
For an algebraic manifold $Z$ over $\C$ we denote the corresponding complex manifold by $Z_\an$, and for a morphism $f:Z\to Z'$ of algebraic manifolds we denote the corresponding holomorphic map by $f_\an:Z_\an\to Z'_\an$.
For an algebraic manifold $Z$ and an $\O_Z$-module $\shf$ we set $\shf_\an=\O_{Z_\an}\tens[\O_Z]\shf$.

In the correspondence \eqref{eq:corr}, let us consider also a functor in the
derived category $\BDC(\C_{\cdot})$ of sheaves of $\C$-vector spaces, going in
the opposite direction:
$$
r:\BDC(\C_{Y_\an})\to \BDC(\C_{X_\an}), \qquad r(F)=\roim{g_\an}\opb{f_\an}(F).
$$

For example, let $D$ be a Zariski locally closed subset of $Y_\an$ and take $F=\C_D$ (the
constant sheaf with fiber $\C$ on $D$ and zero on $Y_\an\setminus D$): then, for any
$x\in X$ one has
\begin{equation}
\label{eq:rfiber}
r(\C_D)_x \simeq \rsect_c(S_{D,x};\C_{S_{D,x}}), \qquad 
S_{D,x}=\{y\in D:(x,y)\in S\}.
\end{equation}

One has the following  ``adjunction formulas'' (see \cite{D'A-S}).
\begin{proposition}
\label{pr:adjunction} Let $\shl$ be a locally free $\O_X$-module
and let $F\in\BDC(\C_{Y_\an})$. Then, setting
$l=\dim Y-\dim S$ and $m=\dim S +\dim Y -2\dim X$, we have
\begin{eqnarray}
\label{eq:af1}
\rsect(X_\an;r(F)\tens\shl_\an^*) &\simeq&
\RHom[\D_{Y_\an}](R(\D\shl)_\an,F\tens\O_{Y_\an}) [l], \\
\label{eq:af2}
\RHom(r(F),\shl_\an^*) &\simeq&
\RHom[\D_{Y_\an}](R(\D\shl)_\an,\rhom(F,\O_{Y_\an}))
[m].
\end{eqnarray}
\end{proposition}

Once the calculation of $R(\D\shl)$ has been performed, these
formulas will give different applications by computing $r(F)$ for different
choices of the sheaf $F$ (a problem of geometric nature).

\subsection{Quasi-equivariant $\D$-modules}
Let us recall the definition of (quasi-)equivariant
$\D$-modules (we refer to Kashiwara~\cite{Kas}). 

Let $G$ be an algebraic group over $\C$, and let $\gothg$ be its Lie algebra.
We denote the enveloping algebra of $\gothg$ by $\shu(\gothg)$.
Let $Z$ be a $G$-manifold, i.e.\ an algebraic manifold endowed with an action of $G$. 
Let us denote by $\mu:G\times Z\to Z$ the action $\mu(g,z)=gz$ and by $p:G\times Z\to Z$ the projection $p(g,z)=z$.
Moreover, define the morphisms
$q_j:G\times G\times Z\to G\times Z$ $(j=1,2,3)$ by $q_1(g_1,g_2,z)=(g_1,g_2z)$,
$q_2(g_1,g_2,z)=(g_1g_2,z)$ and $q_3(g_1,g_2,z)=(g_2,z)$, and observe that $\mu\circ q_1=\mu\circ q_2$, $p\circ q_2=p\circ q_3$ and $\mu\circ q_3=p\circ q_1$. 

A $G$-equivariant $\O_Z$-module is an $\O_Z$-module $\shm$ endowed with an $\O_{G\times Z}$-linear isomorphism 
$\beta:{\mu}^*\shm\to{p}^*\shm$ 
such that the following diagram commutes: 
$$
\begin{array}{ccc}
{q_2}^*{\mu}^*\shm
&
\mathop{\hbox to 170pt{\rightarrowfill}}
\limits^{{q_2}^*\beta}
&
{q_2}^*{p}^*\shm
\\
\downeq
&
&
\downeq
\\
{q_1}^*{\mu}^*\shm
&
\maprightu{{q_1}^*\beta}
{q_1}^*{p}^*\shm
\cong
{q_3}^*{\mu}^*\shm
\maprightu{{q_3}^*\beta}
&
{q_3}^*{p}^*\shm.
\end{array}
$$
For a $G$-equivariant $\O_Z$-module $\shm$ we have a canonical Lie algebra homomorphism
$\rho_{\shm}:\gothg\to\End{\C}(\shm)$.

Let $\O_G\etens\D_Z$
denote the subalgebra $\O_{G\times Z}\tens[\opb{p}\O_Z]
\opb{p}\D_Z$ of $\D_{G\times Z}$. 
A $\D_Z$-module $\shm$ is called $G$-equivariant
(resp.\ quasi-$G$-equivariant) if it is endowed with a $G$-equivariant $\O_Z$-module structure such that the isomorphism
$\beta:{\mu}^*\shm\to{p}^*\shm$  is 
$\D_{G\times Z}$-linear (resp.\ $\O_G\etens\D_Z$-linear).
Note that for a morphism $f:Z\to Z'$ of algebraic manifolds and a $\D_{Z'}$-module $\shm$ the $\D_Z$-module
$H^0(\dopb{f}\shm)$ is naturally isomorphic to $f^*\shm$ as an $\O_Z$-module.

For example for a $G$-equivariant $\O_Z$-module $\shf$ the $\D_Z$-module $\D_Z\tens[\O_Z]\shf$ is endowed with a natural quasi-$G$-equivariant $\D_Z$-module structure.

We denote by $\Mod_G(\D_Z)$ the category of quasi-$G$-equivariant $\D_Z$-modules, and by $\BDC_G(\D_Z)$ the derived category of $\D_Z$-modules with bounded quasi-$G$-equivariant cohomology (see Kashiwara-Schmid~\cite{K-Sc}).

Let $\shm$ be a quasi-$G$-equivariant $\D_Z$-module.
The canonical Lie algebra homomorphism $\gothg\to \D_Z$ induces a Lie algebra homomorphism
$\kappa_{\shm}:\gothg\to\End{\C}(\shm)$.
Set $\gamma_{\shm}=\rho_{\shm}-\kappa_{\shm}$.
\begin{proposition}
[Kashiwara~\cite{Kas}]
\begin{itemize}
\item[\rm(i)]
We have $\gamma_{\shm}(a)\in\End{\D_Z}(\shm)$ for any $a\in\gothg$.
\item[\rm(ii)]
The linear map $\gamma_{\shm}:\gothg\to\End{\D_Z}(\shm)$ is a Lie algebra homomorphism.
\item[\rm(iii)]
We have $\gamma_{\shm}=0$ if and only if  $\shm$ is $G$-equivariant.
\end{itemize}
\end{proposition}
We also denote by 
\begin{equation}
\gamma_{\shm}:\shu(\gothg)\to\End{\D_Z}(\shm)
\end{equation} 
the corresponding algebra homomorphism.

Fix $x\in Z$ and set $H=\{g\in G\,:\,gx=x\}$.
For a $G$-equivariant $\O_Z$-module $\shm$, the fiber 
$$\shm(x)=\C\tens[\O_{Z,x}]\shm_x$$ 
of $\shm$ at $x$ is endowed with a natural $H$-module structure.
If $\shm$ is a quasi-$G$-equivariant $\D_Z$-module, then $\shm(x)$ is also endowed with a $\gothg$-module structure induced from the $\O_Z$-linear action $\gamma_\shm$.
For $M=\shm(x)$ we have the following.
\begin{itemize}
\item[(a)]
the action of the Lie algebra of $H$ on $M$ given by differentiating the $H$-module structure coincides with the restriction of the action of $\gothg$,
\item[(b)]
$hum=(\Ad(h)u)hm$ for any $h\in H$, $u\in\gothg$, $m\in M$.
\end{itemize}
Here $\Ad$ denotes the adjoint action.
A vector space $M$ equipped with structures of an $H$-modules and a $\gothg$-module is called a $(\gothg,H)$-module if it satisfies the conditions (a) and (b) above.

The following result plays a crucial role in the rest of this paper.
\begin{proposition}
\label{pr:G/H}
Assume that $Z=G/H$, where $H$ is a closed subgroup of $G$, and set $x=eH\in Z$.
\begin{itemize}
\item[\rm(i)]
The category of $G$-equivariant $\O_Z$-modules is  equivalent to the category of $H$-modules via the correspondence $\shm\mapsto\shm(x)$.
\item[\rm(ii)]
The category of quasi-$G$-equivariant $\D_Z$-modules is  equivalent to the category of $(\gothg,H)$-modules via the correspondence $\shm\mapsto\shm(x)$.
\end{itemize}
\end{proposition}
The statement (i) is well-known (see \cite{Mumford}), and (ii) is due to Kashiwara~\cite{Kas}.

\section{Radon transforms for generalized flag manifolds}
\subsection{Quasi-equivariant $\D$-modules on generalized flag manifolds}
\label{susec:qeD on gfm}
Let $G$ be a connected reductive algebraic group over $\C$,
and $\gothg$ the Lie algebra of $G$. The group $G$ acts on $\gothg$ by the
adjoint action $\Ad$. 
Let $\gothh$ be a Cartan subalgebra of $\gothg$, $\Delta$ the root
system in $\gothh^*$, $\{\alpha_i:i\in I_0\}$ a set of simple roots, 
$\Delta^+$ the set of positive roots, $\Delta^-$ the set of negative roots, $\gothh^*_{\Z}=\Hom(H,\C^{\times})\subset\gothh^*$ the weight lattice, and $W$ the Weyl group.
For $\alpha\in\Delta$ we denote by $\gothg_\alpha$ the corresponding root space and by $\alpha^\vee\in\gothh$ the corresponding coroot.
For $i\in I_0$ we denote by
$s_i\in W$ the reflection corresponding to $i$.
For $w\in W$ we set $\ell(w)=\sharp(w\Delta^-\cap\Delta^+)$.
Set $\rho=\frac12\sum_{\alpha\in\Delta^+}\alpha$, and define a (shifted) affine action of $W$ on $\gothh^*$ by 
\begin{equation}
w\circ\lambda=w(\lambda+\rho)-\rho. 
\end{equation}

For $I\subset I_0$, we set  
{\small
\begin{eqnarray*}
\Delta_I &=&\Delta\cap\sum_{i\in I}\Z\alpha_i, \qquad
\Delta^+_I=\Delta_I\cap\Delta^+, \qquad  W_I=\langle s_i:i\in I\
\rangle \subset W \\
\gothl_I &=&
\gothh\oplus\left(\bigoplus\nolimits_{\alpha\in\Delta_I}\gothg_{\alpha}
\right), \qquad 
\gothn_I=\bigoplus\nolimits_{\alpha\in\Delta^{+}\setminus
\Delta_I}\gothg_{\alpha}, \qquad 
\gothp_I=\gothl_I\oplus\gothn_I,\\
(\gothh^*_{\Z})_I &=& \{\lambda\in\gothh^*_{\Z}:\lambda(\alpha^\vee_i)\geqq 0 
\mbox{ for any $i\in I$}\}, \\
(\gothh^*_{\Z})_I^0 &=& \{\lambda\in\gothh^*_{\Z}:\lambda(\alpha^\vee_i)= 0 
\mbox{ for any $i\in I$}\}\subset (\gothh^*_{\Z})_I,\\
\rho_I&=&(\sum_{\alpha\in\Delta^+\setminus\Delta_I}\alpha)/2.
\end{eqnarray*}
}We denote by $w_I$ the longest element of $W_I$.
It is an element of $W_I$ characterized by $w_I(\Delta_I^-)=\Delta_I^+$.
Let $L_I$, $N_I$ and ${P_I}$ be the subgroups of $G$ corresponding
to $\gothl_I$, $\gothn_I$ and $\gothp_I$.

For $\lambda\in (\gothh^*_{\Z})_I$ let $V_I(\lambda)$ be the irreducible ${L_I}$-module with highest weight $\lambda$. 
We regard $V_I(\lam)$ as a $P_I$-module with the trivial action of $N_I$, and define the generalized Verma module with highest weight $\lambda$ by 
\begin{equation}
M_I(\lambda)=\shu(\gothg)\tens[\shu(\sgothp_I)]V_I(\lambda).
\end{equation}
Let $L(\lambda)$ be the unique irreducible quotient of $M_I(\lambda)$ (note that $L(\lambda)$ does not depend on the choice of $I$ such that $\lambda\in(\gothh^*_{\Z})_I$).
Then  any irreducible $P_I$-module is isomorphic to $V_I(\lam)$ for some $\lambda\in(\gothh^*_{\Z})_I$, and we have $\dim V_I(\lambda)=1$ if and only if $\lambda\in(\gothh^*_{\Z})_I^0$.
Moreover, any irreducible $(\gothg,P_I)$-module is isomorphic to $L(\lambda)$ for some $\lambda\in(\gothh^*_{\Z})_I$.

Let 
$$
X_I=G/{P_I}
$$
be the generalized flag manifold associated to $I$.

By the category equivalence given in Proposition~\ref{pr:G/H} isomorphism classes of $G$-equivariant $\O_{X_I}$-modules (resp.\ quasi-$G$-equivariant $\D_{X_I}$-modules) are in one-to-one correspondence with isomorphism classes of $P_I$-modules (resp.\ $(\gothg,P_I)$-modules).
For $\lambda\in(\gothh^*_{\Z})_I$ we denote by $\O_{X_I}(\lambda)$ the $G$-equivariant $\O_{X_I}$-module corresponding to the irreducible $P_I$-module $V_I(\lambda)$.
We see easily the following.
\begin{lemma}
\label{lem:qe}
Let $\lambda\in(\gothh^*_{\Z})_I$.
The quasi-$G$-equivariant $\D_{X_I}$-module corresponding to the $(\gothg,P_I)$-module $M_I(\lambda)$ is isomorphic to $\D\O_{X_I}(\lambda)=\D_{X_I}\tens[\O_{X_I}]\O_{X_I}(\lambda)$.
\end{lemma}

We need the following relative version of the Borel-Weil-Bott theorem later (see Bott~\cite{Bott}).
\begin{proposition}
\label{pr:BWW}
Let $I\subset J\subset I_0$ and let $\pi:X_I\to X_J$ be the canonical projection.
For $\lambda\in(\gothh^*_\Z)_I$ we have the following.
\begin{itemize}
\item[\rm(i)]
If there exists some $\alpha\in\Delta_J$ satisfying $(\lambda+\rho-2\rho_I)(\alpha^\vee)=0$, then we have $R\pi_*(\O_{X_I}(\lambda))=0$.
\item[\rm(ii)]
Assume that $(\lambda+\rho-2\rho_I)(\alpha^\vee)\ne0$ for any $\alpha\in\Delta_J$.
Take $w\in W_J$ satisfying $(w(\lambda+\rho-2\rho_I))(\alpha^\vee)>0$ for any $\alpha\in\Delta^+_J$.
Then we have 
$$
R\pi_*(\O_{X_I}(\lambda))=
\O_{X_J}(w(\lambda+\rho-2\rho_I)-(\rho-2\rho_J))
[-(\ell(w_Jw)-\ell(w_I))].
$$
\end{itemize}
\end{proposition}

Let $I,J\subset I_0$ with $I\neq J$. The diagonal action of
$G$ on $X_I\times {X_J}$ has a finite number of orbits, and the only closed
one $G(eP_I,eP_J)$ is identified with ${X_{I\cap J}}=G/({P_I}\cap {P_J})$. 
In the rest of this paper we shall consider the correspondence \eqref{eq:corr} for $X=X_I$, $Y=X_J$ and $S=X_{I\cap J}$:
\begin{equation}
\label{eq:homcorr}
X_I\overset{f}{\longleftarrow} X_{I\cap J}\overset{g}{\longrightarrow}X_J
\end{equation}
and the Radon transform $R(\D\O_{X_I}(\lambda))$ for $\lambda\in(\gothh^*_\Z)_I$.
Since  $f$ and $g$ are morphisms of $G$-manifolds,
the functor \eqref{eq:Radon} induces a functor
\begin{equation}
\label{eq:homRadon}
R:\BDC_G(\D_{X_I}) \to \BDC_G(\D_{X_J}).
\end{equation}
Note that we have
\begin{equation}
\label{eq:gamma}
\Omega_g\simeq\O_{X_{I\cap J}}(\gamma_{I,J})\qquad
\mbox{for $\gamma_{I,J}=\sum_{\alpha\in\Delta_J^+\setminus\Delta_I}\alpha$.}
\end{equation}

\subsection{Radon transforms of quasi-equivariant $\D$-modules}
\label{subsection:RTQED}
Let $\lambda\in(\gothh^*_{\Z})_I$.
We describe our method to analyze $R(\D\O_{X_I}(\lambda))=\doim{g}\dopb{f}(\D\O_{X_I}(\lambda))$.

By 
$$
(\dopb{f}(\D\O_{X_I}(\lambda)))(e(P_I\cap P_J))
\simeq\D\O_{X_I}(\lambda)(eP_I)\simeq M_I(\lam)
$$
the quasi-$G$-equivariant $\D_{X_{I\cap J}}$-module $\dopb{f}(\D\O_{X_I}(\lambda))$ corresponds to the $(\gothg,P_I\cap P_J)$-module $M_I(\lam)=\shu(\gothg)\tens[\shu(\sgothp_I)]V_I(\lambda)$ under the category equivalence given in Proposition~\ref{pr:G/H}.

Set 
\begin{align}
\Gamma&=
\{x\in W_I\,:\,\mbox{$x$ is the shortest element of $W_{I\cap J}x$}\},
\\
\Gamma_k&=\{x\in \Gamma\,:\,\ell(x)=k\}.
\end{align}
It is well-known that an element $x\in W_I$ belongs to $\Gamma$ if and only if $x^{-1}\Delta^+_{I\cap J}\subset\Delta^+_I$.
This condition is also equivalent to 
\begin{equation}
\label{eq:cond:xi1}
\mbox{
$(x(\lambda+\rho))(\alpha^\vee)>0$ for any $\alpha\in\Delta^+_{I\cap J}$.
}
\end{equation}
In particular, we have $x\circ\lambda\in(\gothh_\Z^*)_{I\cap J}$ for $x\in\Gamma$.

By Lepowsky~\cite{Lep} and Rocha-Caridi~\cite{Rocha-Caridi} we have the following resolution of the finite dimensional $\gothl_I$-module $V_I(\lambda)$:
\begin{equation}
\label{eq:BGG1}
0\to N_n\to N_{n-1}\to\cdots \to N_1\to N_0\to V_I(\lambda)\to0
\end{equation}
with $n=\dim\gothl_I/\gothl_I\cap\gothp_J$ and 
$$
N_k=\bigoplus_{x\in\Gamma_k}\shu(\gothl_I)\tens[\shu(\sgothl_I\cap\sgothp_J)]V_{I\cap J}(x\circ\lambda).
$$
By the Poincar\'e-Birkhoff-Witt theorem we have the isomorphism 
$$
\shu(\gothl_I)\tens[\shu(\sgothl_I\cap\sgothp_J)]V_{I\cap J}(x\circ\lambda)
\simeq
\shu(\gothp_I)\tens[\shu(\sgothp_{I\cap J})]V_{I\cap J}(x\circ\lambda)
$$
of $\shu(\gothl_I)$-modules, where $\gothn_{I\cap J}$ acts trivially on $V_{I\cap J}(x\circ\lambda)$.
Moreover, the action of $\gothn_I$ on $\shu(\gothp_I)\tens[\shu(\sgothp_{I\cap J})]V_{I\cap J}(x\circ\lambda)$ is trivial.
Indeed, by $[\gothp_I,\gothn_I]\subset\gothn_I$ we have $\gothn_I\shu(\gothp_I)=\shu(\gothp_I)\gothn_I$, and hence
$$
\gothn_I
(\shu(\gothp_I)\tens[\shu(\sgothp_{I\cap J})]V_{I\cap J}(x\circ\lambda))
\subset
\shu(\gothp_I)\gothn_I\tens V_{I\cap J}(x\circ\lambda)
\subset
\shu(\gothp_I)\tens \gothn_IV_{I\cap J}(x\circ\lambda)
=0
$$
by $\gothn_I\subset\gothn_{I\cap J}$.
Thus we obtain the following resolution of the finite dimensional $\gothp_I$-module $V_I(\lambda)$ (with trivial action of $\gothn_I$):
\begin{equation}
\label{eq:BGG2}
0\to N'_n\to N'_{n-1}\to\cdots \to N'_1\to N'_0\to V_I(\lambda)\to0
\end{equation}
with 
$$
N'_k=\bigoplus_{x\in\Gamma_k}\shu(\gothp_I)\tens[\shu(\sgothp_{I\cap J})]V_{I\cap J}(x\circ\lambda).
$$
By tensoring $\shu(\gothg)$ to \eqref{eq:BGG2} over $\shu(\gothp_I)$ we obtain the following resolution of the $(\gothg, P_{I\cap J})$-module $M_I(\lambda)$:
\begin{equation}
\label{eq:BGG3}
0\to \tilde{N}_n\to \tilde{N}_{n-1}\to\cdots \to \tilde{N}_1\to \tilde{N}_0\to M_I(\lambda)\to0
\end{equation}
with 
$$
\tilde{N}_k=\bigoplus_{x\in\Gamma_k}M_{I\cap J}(x\circ\lambda).
$$
Since the quasi-$G$-equivariant $\D_{X_{I\cap J}}$-module corresponding to the $(\gothg, P_{I\cap J})$-module $M_{I\cap J}(x\circ\lambda)$ is $\D\O_{X_{I\cap J}}(x\circ\lambda)$, we have obtained the following resolution of the quasi-$G$-equivariant $\D_{X_{I\cap J}}$-module $\dopb{f}(\D\O_{X_I}(\lambda))$:
\begin{equation}
\label{eq:BGG4}
0\to \shn_n\to \shn_{n-1}\to\cdots \to \shn_1\to \shn_0\to\dopb{f}(\D\O_{X_I}(\lambda))\to0
\end{equation}
with 
\begin{equation}
\label{eq:BGG4+}
\shn_k=\bigoplus_{x\in\Gamma_k}\D\O_{X_{I\cap J}}(x\circ\lambda).
\end{equation}

Our next task is to investigate on $\doim{g}(\D\O_{X_{I\cap J}}(x\circ\lambda))$ for $x\in \Gamma$.
We first remark that 
\begin{equation}
\label{eq:BGG5}
\doim{g}(\D\O_{X_{I\cap J}}(x\circ\lambda))
=\D_{X_J}\tens[\O_{X_J}]Rg_*(\O_{X_{I\cap J}}(x\circ\lambda+\gamma_{I,J})).
\end{equation}
Indeed, by \eqref{eq:gamma} we have
\begin{align*}
\doim{g}(\D\O_{X_{I\cap J}}(x\circ\lambda))
&=
Rg_*(\D_{X_J\from X_{I\cap J}}\ltens[\D_{X_{I\cap J}}]\D_{X_{I\cap J}}\ltens[\O_{X_{I\cap J}}]\O_{X_{I\cap J}}(x\circ\lambda))\\
&=
Rg_*(\D_{X_J\from X_{I\cap J}}\ltens[\O_{X_{I\cap J}}]\O_{X_{I\cap J}}(x\circ\lambda))\\
&=
Rg_*(g^{-1}\D_{X_J}\tens[g^{-1}\O_{X_J}]\Omega_g\tens[\O_{X_{I\cap J}}]
\O_{X_{I\cap J}}(x\circ\lambda))\\
&=
\D_{X_J}\tens[\O_{X_J}]Rg_*(\Omega_g\tens[\O_{X_{I\cap J}}]
\O_{X_{I\cap J}}(x\circ\lambda))\\
&=
\D_{X_J}\tens[\O_{X_J}]Rg_*(\O_{X_{I\cap J}}(x\circ\lambda+\gamma_{I,J})).
\end{align*}

\begin{lemma}
\label{le:BWW2}
Let $\lambda\in(\gothh_\Z^*)_I$ and $x\in\Gamma$.
\begin{itemize}
\item[\rm(i)]
If $(x(\lambda+\rho))(\alpha^\vee)=0$ for some $\alpha\in\Delta_J$, then we have $Rg_*(\O_{X_{I\cap J}}(x\circ\lambda+\gamma_{I,J}))=0$.
\item[\rm(ii)]
Assume that $(x(\lambda+\rho))(\alpha^\vee)\ne0$ for any $\alpha\in\Delta_J$.
Take $y\in W_J$ satisfying $(yx(\lambda+\rho))(\alpha^\vee)>0$ for any $\alpha\in\Delta^+_J$.
Then we have
$$
Rg_*(\O_{X_{I\cap J}}(x\circ\lambda+\gamma_{I,J}))=
\O_{X_J}((yx)\circ\lambda)[-(\ell(w_Jy)-\ell(w_{I\cap J}))].
$$
\end{itemize}
\end{lemma}
\begin{proof}
Since $\Delta^+\setminus\Delta_J$ is stable under the action of $W_J$, we have $y\rho_{J}=\rho_{J}$ for any $y\in W_J$.
In particular, 
$$
\rho_{J}=s_\alpha(\rho_{J})=\rho_J-\rho_J(\alpha^\vee)\alpha
$$
for any $\alpha\in\Delta_J$, and hence $\rho_J(\alpha^\vee)=0$ for any $\alpha\in\Delta_J$.

By the definition we have
$$
x\circ\lambda+\gamma_{I,J}+\rho-2\rho_{I\cap J}
=x(\lambda+\rho)+\gamma_{I,J}-2\rho_{I\cap J}
=x(\lambda+\rho)-2\rho_{J},
$$
and 
$$
y(x(\lambda+\rho)-2\rho_{J})-(\rho-2\rho_J)
=yx(\lambda+\rho)-2\rho_{J}-(\rho-2\rho_J)
=(yx)\circ\lambda
$$
for any $y\in W_J$.
Hence the assertion follows from Proposition~\ref{pr:BWW}.
\end{proof}
Set 
\begin{align}
\Gamma(\lambda)&=\{x\in\Gamma\,:\,
\mbox{$(x(\lambda+\rho))(\alpha^\vee)\ne0$ for any $\alpha\in\Delta_J$}
\},\\
\Gamma_k(\lambda)&=\{x\in \Gamma(\lambda)\,:\,\ell(x)=k\}.
\end{align}
and for $x\in\Gamma(\lambda)$ denote by $y_x$ the element of $W_J$ satisfying $(y_xx(\lambda+\rho))(\alpha^\vee)>0$ for any $\alpha\in\Delta^+_J$.
Set 
\begin{equation}
m(x)=\ell(w_Jy_x)-\ell(w_{I\cap J})\qquad
\mbox{for $x\in\Gamma(\lambda)$.}
\end{equation}
\begin{lemma}
\label{lem:lm}
For $\lambda\in(\gothh_\Z^*)_I$ and $x\in\Gamma(\lambda)$ we have
\begin{align}
\ell(x)&=\sharp\{\alpha\in\Delta^+_I\setminus\Delta_J\,:\,
(x(\lambda+\rho))(\alpha^\vee)<0\},\\
m(x)&=\sharp\{\alpha\in\Delta_J^+\setminus\Delta_I\,:\,
(x(\lambda+\rho))(\alpha^\vee)>0\}.
\end{align}
\end{lemma}
\begin{proof}
We have
\begin{align*}
\ell(x)&=
\sharp(x^{-1}\Delta_I^-\cap \Delta_I^+)\\
&=\sharp\{\alpha\in\Delta^+_I\,:\,
(x(\lambda+\rho))(\alpha^\vee)<0\}\\
&=\sharp\{\alpha\in\Delta^+_I\setminus\Delta_J\,:\,
(x(\lambda+\rho))(\alpha^\vee)<0\},
\end{align*}
and
\begin{align*}
m(x)&=
\ell(w_J)-\ell(y_x)-\ell(w_{I\cap J})\\
&=\sharp(\Delta_J^+\setminus\Delta_I)-\sharp(y_x^{-1}\Delta_J^-\cap \Delta_J^+)\\
&=\sharp(\Delta_J^+\setminus\Delta_I)
-\sharp\{\alpha\in\Delta^+_J\,:\,
(x(\lambda+\rho))(\alpha^\vee)<0\}\\
&=\sharp(\Delta_J^+\setminus\Delta_I)
-\sharp\{\alpha\in\Delta^+_J\setminus\Delta_I\,:\,
(x(\lambda+\rho))(\alpha^\vee)<0\}\\
&=\sharp\{\alpha\in\Delta^+_J\setminus\Delta_I\,:\,
(x(\lambda+\rho))(\alpha^\vee)>0\}
\end{align*}
by \eqref{eq:cond:xi1}.
\end{proof}
\begin{proposition}
\label{prop:Mk}
Let $\lambda\in(\gothh_\Z^*)_I$.
Then there exists a family $\{\shm(k)^\bullet\}_{k\geqq0}$ of objects of $\BDC_G(\D_{X_J})$ satisfying the following conditions.
\begin{itemize}
\item[\rm(i)]
$\shm(0)^\bullet\simeq R(\D\O_{X_I}(\lambda))$.
\item[\rm(ii)]
$\shm(k)^\bullet=0$ for $k>\dim\gothl_I/\gothl_I\cap\gothp_J$.
\item[\rm(iii)]
We have a distinguished triangle
\[
\shc(k)^\bullet\to
\shm(k)^\bullet\to
\shm(k+1)^\bullet\stackrel{+1}{\longrightarrow}
\]
where
\[
\shc(k)^\bullet=
\bigoplus_{x\in\Gamma_k(\lambda)}
\D\O_{X_J}((y_xx)\circ\lambda)[\ell(x)-m(x)].
\]
\end{itemize}
\end{proposition}
\begin{proof}
For $0\leqq k\leqq \dim\gothl_I/\gothl_I\cap\gothp_J$ define an object $\shn(k)^\bullet$ of $\BDC_G(\D_{X_{I\cap J}})$ by
$$
\shn(k)^\bullet=
[\cdots\to0\to\shn_n\to \shn_{n-1}\to\cdots \to \shn_k\to0\cdots],
$$
where $\shn_j$ has degree $-j$
(see \eqref{eq:BGG4} and \eqref{eq:BGG4+} for the notation).
For $k>\dim\gothl_I/\gothl_I\cap\gothp_J$ we set $\shn(k)^\bullet=0$.
By $\shn(0)^\bullet\simeq\dopb{f}(\D\O_{X_I}(\lambda))$ we have
$\doim{g}\shn(0)^\bullet\simeq R(\D\O_{X_I}(\lambda))$.
Set $\shm(k)^\bullet=\doim{g}\shn(k)^\bullet$.
Then the statements (i) and (ii) are obvious.
Let us show (iii).
Applying $\doim{g}$ to the distinguished triangle
$$
\shn_k[k]\to
\shn(k)^\bullet\to
\shn(k+1)^\bullet\stackrel{+1}{\longrightarrow}
$$
we obtain a distinguished triangle
$$
\doim{g}\shn_k[k]\to
\shm(k)^\bullet\to
\shm(k+1)^\bullet\stackrel{+1}{\longrightarrow}.
$$
By \eqref{eq:BGG4+}, \eqref{eq:BGG5} and Lemma~\ref{le:BWW2} we have
$$
\doim{g}\shn_k=
\bigoplus_{x\in\Gamma_k(\lambda)}
\D\O_{X_J}((y_xx)\circ\lambda)[-m(x)].
$$
The statement (iii) is proved.
\end{proof}
\begin{theorem}
\label{thm:main}
Let $\lambda\in(\gothh_\Z^*)_I$.
\begin{itemize}
\item[\rm(i)]
We have
$$
\sum_p(-1)^p[H^p(R(\D\O_{X_I}(\lambda)))]
=\sum_{x\in\Gamma(\lambda)}(-1)^{\ell(x)-m(x)}[\D\O_{X_J}((y_xx)\circ\lambda)]
$$
in the Grothendieck group of the category of quasi-$G$-equivariant $\D_{X_J}$-modules.
\item[\rm(ii)]
If $\Gamma(\lambda)=\emptyset$, then $R(\D\O_{X_I}(\lambda))=0$.
\item[\rm(iii)]
If $\Gamma(\lambda)$ consists of a single element $x$, then 
$$
R(\D\O_{X_I}(\lambda))=\D\O_{X_J}((y_xx)\circ\lambda)[\ell(x)-m(x)].
$$
\item[\rm(iv)]
If $\ell(x)\geqq m(x)$ for any $x\in\Gamma(\lambda)$, then we have $H^p(R(\D\O_{X_I}(\lambda)))=0$ unless $p=0$.
\item[\rm(v)]
If $(\lambda+\rho)(\alpha^\vee)<0$ for any $\alpha\in\Delta^+_J\setminus\Delta_I$, then there exists a canonical morphism
$$
\Phi:\D\O_{X_J}((w_Jw_{I\cap J})\circ\lambda)\to H^0(R(\D\O_{X_I}(\lambda))).
$$
Moreover, $\Phi$ is an epimorphism if $\ell(x)>m(x)$ for any $x\in\Gamma(\lambda)\setminus\{e\}$.
\end{itemize}
\end{theorem}
\begin{proof}
The statements (i), (ii), (iii) are obvious from Proposition~\ref{prop:Mk}.
The statement (iv) follows from Proposition~\ref{prop:Mk} and Lemma~\ref{lem:lower vanishing}.
Assume that $\lambda$ satisfies the assumption in (v).
Then we have $e\in\Gamma(\lambda)$ and $y_e=w_Jw_{I\cap J}$.
Hence (v) follows from Proposition~\ref{prop:Mk}.
\end{proof}
\begin{lemma}
\label{lem:WJWI}
\begin{itemize}
\item[\rm(i)]
The map $W_J\times\Gamma\to W_JW_I\,((y,x)\mapsto yx)$ is bijective.
\item[\rm(ii)]
For $\lambda\in(\gothh_\Z^*)_I$ we have
$$
\{y_xx\,:\,x\in\Gamma(\lambda)\}
=\{w\in W_JW_I\,:\,
\mbox{$(w(\lambda+\rho))(\alpha^\vee)>0$ for any $\alpha\in\Delta^+_J$}
\}
$$
and we have 
$$
\ell(x)-m(x)=\ell(y_x)+\ell(x)-\sharp(\Delta_J^+\setminus\Delta_I)=\ell(y_xx)-\sharp(\Delta_J^+\setminus\Delta_I).
$$
\end{itemize}
\end{lemma}
\begin{proof}
(i) is a consequence of the definition of $\Gamma$, and the first statement in (ii) follows from (i) and the definition of $y_x$.
By
$$
\ell(x)-m(x)=\ell(x)-(\ell(w_J)-\ell(y_x)-\ell(w_{I\cap J}))
=\ell(x)+\ell(y_x)-\sharp(\Delta_J^+\setminus\Delta_I)
$$
we have only to show $\ell(y_xx)=\ell(x)+\ell(y_x)$ for $x\in\Gamma(\lambda)$.
We have 
$$
x\Delta^+\cap\Delta^-=x\Delta_I^+\cap\Delta_I^-\subset\Delta_I^-\setminus\Delta_{I\cap J}\subset\Delta^-\setminus\Delta_J
$$
by $x\in W_I$ and $x^{-1}\Delta_{I\cap J}^+\subset\Delta^+_I$.
Since $w\in W_J$, we obtain $y_x(x\Delta^+\cap\Delta^-)\subset\Delta^-$.
Hence
\begin{align*}
\ell(y_xx)&=
\sharp(y_xx\Delta^-\cap\Delta^+)\\
&=\sharp(y_x(x\Delta^-\cap\Delta^+)\cap\Delta^+)
+\sharp(y_x(x\Delta^-\cap\Delta^-)\cap\Delta^+)\\
&=\sharp(y_x(x\Delta^-\cap\Delta^+)\cap\Delta^+)
+\sharp(y_x\Delta^-\cap\Delta^+)
-\sharp(y_x(x\Delta^+\cap\Delta^-)\cap\Delta^+)\\
&=\ell(x)+\ell(y_x).
\end{align*}
\end{proof}
For $\lambda\in(\gothh_\Z^*)_I$ we set
\begin{equation}
\Xi(\lambda)
=\{w\in W_JW_I\,:\,
\mbox{$(w(\lambda+\rho))(\alpha^\vee)>0$ for any $\alpha\in\Delta^+_J$}\}.
\end{equation}
Using Lemma~\ref{lem:WJWI} above we can reformulate Theorem~\ref{thm:main} as follows.
\begin{theorem}
\label{thm:main+}
Let $\lambda\in(\gothh_\Z^*)_I$.
\begin{itemize}
\item[\rm(i)]
We have
$$
\sum_p(-1)^p[H^p(R(\D\O_{X_I}(\lambda)))]
=(-1)^{\sharp(\Delta^+_J\setminus\Delta_I)}\sum_{w\in\Xi(\lambda)}(-1)^{\ell(w)}[\D\O_{X_J}(w\circ\lambda)]
$$
in the Grothendieck group of the category of quasi-$G$-equivariant $\D_{X_J}$-modules.
\item[\rm(ii)]
If $\Xi(\lambda)=\emptyset$, then $R(\D\O_{X_I}(\lambda))=0$.
\item[\rm(iii)]
If $\Xi(\lambda)$ consists of a single element $w$, then 
$$
R(\D\O_{X_I}(\lambda))=\D\O_{X_J}(w\circ\lambda)[\ell(w)-\sharp(\Delta^+_J\setminus\Delta_I)].
$$
\item[\rm(iv)]
If $\ell(w)\geqq\sharp(\Delta^+_J\setminus\Delta_I)$ for any $w\in\Xi(\lambda)$, then we have $H^p(R(\D\O_{X_I}(\lambda)))=0$ unless $p=0$.
\item[\rm(v)]
If $(\lambda+\rho)(\alpha^\vee)<0$ for any $\alpha\in\Delta^+_J\setminus\Delta_I$, then there exists a canonical morphism
$$
\Phi:\D\O_{X_J}((w_Jw_{I\cap J})\circ\lambda)\to H^0(R(\D\O_{X_I}(\lambda))).
$$
Moreover, $\Phi$ is an epimorphism if $\ell(w)>\sharp(\Delta^+_J\setminus\Delta_I)$ for any $w\in\Xi(\lambda)\setminus\{w_Jw_{I\cap J}\}$.
\end{itemize}
\end{theorem}
\begin{remark}
The following result which is a little weaker than Theorem~\ref{thm:main+}(ii) can be obtained by observing that an integral transform
for $\D$-modules with equivariant kernel preserves the 
infinitesimal character of a quasi-equivariant $\D$-module 
(see e.g.\ \cite{K-Sc}):
\begin{equation}
\mbox{
If $(W\circ\lambda)\cap(\gothh^*_{\Z})_J =\emptyset$,
then  $R(\D\O_{X_I}(\lambda))=0$.}
\end{equation}
An advantage of the argument using the infinitesimal character is that it also works for a broader class of integral transforms in equivariant contexts.

Let us briefly recall this argument (suggested to us by
M.\  Kashiwara). 
Let $Z$ be a $G$-manifold, denote by
$\gothz(\gothg)$ the center of $\shu(\gothg)$ and
set 
$\gothn^{+}=\gothn_{\emptyset}=\bigoplus\nolimits_{
\alpha\in\Delta^{+}}\gothg_{\alpha}$ and $\gothn^-=
\bigoplus\nolimits_{\alpha\in\Delta^{-}}\gothg_{\alpha}$. 
One says that a quasi-$G$-equivariant $\D_Z$-module $\shm$ 
has infinitesimal character $\chi$ (for some $\chi\in
\Hom(\gothz (\gothg),\C)$) if $\gamma_{\shm}(a)$ is the
multiplication by $\chi(a)$ for any $a\in\gothz(\gothg)$.
Define a linear map $\sigma:\gothz(\gothg)\to \shu(\gothh)
\simeq S(\gothh)$ as the composition of the embedding
$\gothz(\gothg)\to\shu(\gothg)$ and the projection
$\shu(\gothg)\to\shu(\gothh)$ with respect to the direct sum
decomposition $\shu(\gothg)=\shu(\gothh)\oplus(\gothn^-
\shu(\gothg) +\shu(\gothg)\gothn^+)$. Then $\sigma$ is an
injective homomorphism of $\C$-algebras. For $\lambda\in
\gothh^*$ define an algebra homomorphism $\chi_{\lambda}:
\gothz(\gothg) \to\C$ by $\chi_{\lambda}(a)=\langle \sigma(a),
\lambda\rangle$. By a result of Harish-Chandra, any algebra
homomorphism from $\gothz(\gothg)$ to $\C$ coincides with
$\chi_\lambda$ for some $\lambda\in\gothh^*$, and for
$\lambda,\mu\in\gothh^*$ one has $\chi_{\lambda}=\chi_{\mu}$
if and only if $\mu\in W\circ\lambda$. By the category
equivalence of Proposition \ref{pr:G/H}, the infinitesimal
characters of quasi-$G$-equivariant $\D_{X_I}$-modules are of
the form $\chi_{\lambda}$ for $\lambda\in(\gothh^*_{\Z})_I$.
Therefore, recalling Harish-Chandra's result, if $(W\circ
\lambda)\cap(\gothh^*_{\Z})_J =\emptyset$, then  $R(\D\O_{X_I}(\lambda))=0$.
\end{remark}

\subsection{Extremal cases}
We characterize the extremal cases (see Definition \ref{de:extrcase})
in the correspondence \eqref{eq:homcorr}. 
We shall only deal with the invertible $\O$-modules.
Given $\lambda\in(\gothh^*_{\Z})_I^0$
and $\mu\in(\gothh^*_{\Z})_J^0$, we write for short $(\lambda,\mu)$ instead of
$(\O_{X_I}(\lambda),\O_{X_J}(\mu))$. 

\begin{proposition}
The pair $(\lambda,\mu)$ is an extremal case if and only if
$\mu=\lambda+\gamma_{I,J}$.
This condition is also equivalent to the following system
\begin{equation}
\label{eq:extrsys}
\left\{
\begin{array}{ll}
\lambda(\alpha_i^\vee)=\mu(\alpha_i^\vee)=0 & (i\in I\cap J), \\
\lambda(\alpha_i^\vee)=0,\,\,\mu(\alpha_i^\vee)=\gamma_{I,J}(\alpha_i^\vee) & (i\in I\setminus J), \\
\lambda(\alpha_i^\vee)=-\gamma_{I,J}(\alpha_i^\vee),\,\,\mu(\alpha_i^\vee)=0 & (i\in J\setminus I), \\
\mu(\alpha_i^\vee)-\lambda(\alpha_i^\vee)=\gamma_{I,J}(\alpha_i^\vee) & (i\in I_0\setminus(I\cup J)).
\end{array}
\right.
\end{equation}
\end{proposition}
\begin{proof}
The first statement is obvious by \eqref{eq:gamma}.
Since $\Delta^+\setminus\Delta_I$ and $\Delta_J$ are stable under the action of $W_I$ and $W_J$ respectively, we have $w(\gamma_{I,J})=\gamma_{I,J}$ for any $w\in W_{I\cap J}=W_I\cap W_J$. 
In particular, we have
$$
\gamma_{I,J}=s_i(\gamma)=\gamma_{I,J}-\gamma_{I,J}(\alpha_i^\vee)\alpha_i
$$
for any $i\in I\cap J$.
Hence we obtain
$$
\gamma_{I,J}(\alpha_i^\vee)=0\quad\mbox{ for any }i\in I\cap J.
$$
Therefore, the relation $\mu=\lambda+\gamma_{I,J}$ is equivalent to the system \eqref{eq:extrsys}.
\end{proof}
By \eqref{eq:extrsys} we have the following
\begin{corollary}
If $\gothg$ is semisimple and $I\cup J=I_0$, there exists a unique extremal case for \eqref{eq:homcorr}.
\end{corollary}
\begin{proposition}
\label{pr:extremal}
If $(\lambda,\mu)$ is an extremal case, then we have
$$
(\lambda+\rho)(\alpha^\vee)
\left\{
\begin{array}{ll}
<0\quad
&\mbox{for any $\alpha\in\Delta_J^+\setminus\Delta_I$},\\
>0\quad
&\mbox{for any $\alpha\in\Delta^+_I$},
\end{array}\right.
$$
and $(w_Jw_{I\cap J})\circ\lambda=\mu$.
In particular, we have $e\in\Gamma(\lambda)$ and $\ell(e)=m(e)=0$.
\end{proposition}
\begin{proof}
Since $\mu$ and $\gamma_{I,J}$ are fixed by the action of $W_J$ and $W_{I\cap J}$ respectively, We have
$$
(w_Jw_{I\cap J})\circ\lambda
=w_Jw_{I\cap J}(\mu-\gamma_{I,J}+\rho)-\rho=\mu-w_J(\gamma_{I,J}-w_{I\cap J}\rho+w_J\rho).
$$
By 
$$
w_{I\cap J}\rho-w_J\rho
=(\rho-w_J\rho)-(\rho-w_{I\cap J}\rho)
=\sum_{\alpha\in\Delta_J^+}\alpha-\sum_{\alpha\in\Delta_{I\cap J}^+}\alpha
=\gamma_{I,J}
$$
we obtain $(w_Jw_{I\cap J})\circ\lambda=\mu$.
Hence by $w_Jw_{I\cap J}(\Delta_J^+\setminus\Delta_I)\subset\Delta_J^-$ and $\mu\in(\gothh^*_\Z)_J^0$, we have
$$
(\lambda+\rho)(\alpha^\vee)
=(w_{I\cap J}w_J(\mu+\rho))(\alpha^\vee)
=(\mu+\rho)(w_Jw_{I\cap J}\alpha^\vee)<0
$$
for any $\alpha\in\Delta_J^+\setminus\Delta_I$.
Moreover, we have $(\lambda+\rho)(\alpha^\vee)>0$ for any $\Delta^+_I$ by \eqref{eq:extrsys}.
\end{proof}

By Proposition \ref{pr:ecnontrmor}, if the pair $(\lambda,\mu)$ is an extremal case we get a nontrivial $\D_{X_J}$-linear morphism
\begin{equation}
\label{eq:Dmorph}
\Phi:\D\O_{X_J}(\mu)\to H^0(R(\D\O_{X_I}(\lambda))).
\end{equation}
\begin{theorem}
\label{thm:extremal case}
Let $(\lambda,\mu)$ be an extremal case.
\begin{itemize}
\item[\rm(i)]
We have $H^p(R(\D\O_{X_I}(\lambda)))=0$ for any $p\ne0$ if and only if $\ell(x)\geqq m(x)$ for any $x\in\Gamma(\lambda)$.
\item[\rm(ii)]
Assume that $H^p(R(\D\O_{X_I}(\lambda)))=0$ for any $p\ne0$.
Then $\Phi$ is an epimorphism if and only if $\ell(x)> m(x)$ for any $x\in\Gamma(\lambda)\setminus\{e\}$.
\item[\rm(iii)]
Assume that $H^p(R(\D\O_{X_I}(\lambda)))=0$ for any $p\ne0$.
Then $\Phi$ is an isomorphism if and only if $\Gamma(\lambda)=\{e\}$.
\end{itemize}
\end{theorem}
We need the following result in order to prove Theorem~\ref{thm:extremal case}.
\begin{lemma}
\label{lem:subquotient}
Let $(\lambda,\mu)$ be an extremal case, and let $x_1, x_2\in\Gamma(\lambda)$.
Set $y_k=y_{x_k}$ for $k=1, 2$.
If $L((y_1x_1)\circ\lambda)$ appears as a subquotient of $M_J((y_2x_2)\circ\lambda)$, then we have $\ell(x_2)-\ell(y_2)\leqq\ell(x_1)-\ell(y_1)$.
\end{lemma}
\begin{proof}
For $\xi\in\gothh_\Z^*$ we set
$$
\Delta_0^+(\xi)=
\{\alpha\in\Delta^+\,:\,
(\xi+\rho)(\alpha^\vee)=0\},\quad
W_0(\xi)=
\{w\in W\,:\,w\circ\xi=\xi\}.
$$
Take $\nu\in W\circ\lambda$ such that $(\nu+\rho)(\alpha^\vee)\geqq0$ for any $\alpha\in\Delta^+$, and let $w\in W$ such that $\lambda=w\circ\nu$.
We can assume that $\ell(w)\leqq\ell(x)$ for any $x\in W$ satisfying $\lambda=x\circ\nu$.
Then $w$ is the (unique) element of $wW_0(\nu)$ with minimal length.

Let us first show:
\begin{equation}
\label{eq:minimal length}
\mbox{$y_kx_kw$ is the element of $y_kx_kwW_0(\nu)$ with minimal length.}
\end{equation}
It is sufficient to show $y_kx_kw\Delta^+_0(\nu)\subset\Delta^+$.
Since $w$ is the element of $wW_0(\nu)$ with minimal length, we have $w\Delta^+_0(\nu)\subset\Delta^+$, and hence $w\Delta^+_0(\nu)=\Delta^+_0(\lambda)$.
By Proposition~\ref{pr:extremal} we have $\Delta^+_0(\lambda)\subset\Delta^+\setminus\Delta_I$.
Hence by $W_I(\Delta^+\setminus\Delta_I)=\Delta^+\setminus\Delta_I$ we have $x_k\Delta^+_0(\lambda)\subset\Delta^+$.
Thus $x_k\Delta^+_0(\lambda)=\Delta^+_0(x_k\circ\lambda)$.
By $x_k\in\Gamma(\lambda)$ we have $\Delta^+_0(x_k\circ\lambda)\subset\Delta^+\setminus\Delta_J$, and hence $y_k\Delta^+_0(x_k\circ\lambda)\subset\Delta^+$.
The statement \eqref{eq:minimal length} is proved.

We next show
\begin{equation}
\label{eq:length}
\ell(y_kx_kw)=\ell(w)+\ell(x_k)-\ell(y_k).
\end{equation}
For any $\alpha\in\Delta_I^+$ we have
$$
(\nu+\rho)(w^{-1}\alpha^\vee)=(\lambda+\rho)(\alpha^\vee)>0,
$$
and hence $w^{-1}\Delta_I^+\subset\Delta^+$ by the choice of $\nu$.
Thus we have 
$$
w^{-1}(x_k^{-1}\Delta^+\cap\Delta^-)=w^{-1}(x_k^{-1}\Delta_I^+\cap\Delta_I^-)
\subset w^{-1}\Delta_I^-\subset\Delta^-.
$$
Hence $\ell(x_kw)=\ell(w)+\ell(x_k)$.
Here, we have used the well-known fact that for $u, v\in W$ we have $\ell(uv)=\ell(u)+\ell(v)$ if and only if $u(v\Delta^+\cap\Delta^-)\subset\Delta^-$.
Similarly, we have
$$
(\nu+\rho)(w^{-1}x_k^{-1}y_k^{-1}\alpha^\vee)
=(y_kx_k(\lambda+\rho))(\alpha^\vee)>0
$$
for any $\alpha\in\Delta_J^+$ by the definition of $y_k$ and hence $w^{-1}x_k^{-1}y_k^{-1}\Delta_J^+\subset\Delta^+$.
Thus we have 
$$
w^{-1}x_k^{-1}y_k^{-1}(y_k\Delta^+\cap\Delta^-)
=w^{-1}x_k^{-1}y_k^{-1}(y_k\Delta_J^+\cap\Delta_J^-)
\subset w^{-1}x_k^{-1}y_k^{-1}\Delta_J^-
\subset \Delta^-.
$$
Hence $\ell(x_kw)=\ell(y_kx_kw)+\ell(y_k)$.
The statement \eqref{eq:length} is proved.

Note that $L((y_1x_1)\circ\lambda)=L((y_1x_1w)\circ\nu)$ and that $M_J((y_2x_2)\circ\lambda)$ is a quotient of the ordinary Verma module $M((y_2x_2w)\circ\nu)=M_\emptyset((y_2x_2w)\circ\nu)$.
Hence by our assumption and by \eqref{eq:minimal length} we obtain $y_1x_1w\geqq y_2x_2w$ with respect to the standard partial order on $W$  by a result of Bernstein-Gelfand-Gelfand~\cite{BGG} concerning the composition factors of Verma modules.
In particular, we have $\ell(y_1x_1w)\geqq \ell(y_2x_2w)$.
Hence we obtain the desired result by \eqref{eq:length}.
\end{proof}
{\sl Proof of Theorem~\ref{thm:extremal case}.}\quad
We shall use the notation in Proposition~\ref{prop:Mk}.

We first show the following.
\begin{equation}
\label{eq:pr of th:1}
\mbox{If $H^r(\shm(k)^\bullet)=0$ for any $k\geqq\ell$, then $H^r(\shc(k)^\bullet)=0$ for any $k\geqq\ell$.}
\end{equation}
Assume that there exists some $k\geqq\ell$ such that $H^r(\shc(k)^\bullet)\ne0$.
Let $k_0$ be the largest such $k$.
Then we have exact sequences
\begin{align}
&H^{r-1}(\shm(k_0+1)^\bullet)\to H^{r}(\shc(k_0)^\bullet)\to 0,
\label{eq:pr of th:2}\\
&H^{r-1}(\shc(k)^\bullet)\to H^{r-1}(\shm(k)^\bullet)\to H^{r-1}(\shm(k+1)^\bullet)\to 0
\qquad(k>k_0).
\label{eq:pr of th:3}
\end{align}
By $H^{r}(\shc(k_0)^\bullet)\ne0$ there exists some $x_1\in\Gamma(\lambda)$ such that $\ell(x_1)-m(x_1)=-r$, $\ell(x_1)=k_0$ and $\D\O_{X_J}((y_{x_1}x_1)\circ\lambda)$ is a direct summand of $H^r(\shc(k_0)^\bullet)$.
On the other hand by \eqref{eq:pr of th:2} and \eqref{eq:pr of th:3}  any irreducible subquotient of $H^r(\shc(k_0)^\bullet)$ is isomorphic to an irreducible subquotient of $H^{r-1}(\shc(k)^\bullet)$ for some $k\geqq k_0+1$.
Moreover, $H^{r-1}(\shc(k)^\bullet)$ is isomorphic to the direct sum of $\D\O_{X_J}((y_{x_2}x_2)\circ\lambda)$ for $x_2\in\Gamma(\lambda)$ such that $\ell(x_2)-m(x_2)=-(r-1)$, $\ell(x_2)=k$.
By the category equivalence given in Proposition~\ref{pr:G/H} we see that there exists some $x_2\in\Gamma(\lambda)$ such that $\ell(x_2)-m(x_2)=-(r-1)$, $\ell(x_2)\geqq k_0+1$, and that $L((y_{x_1}x_1)\circ\lambda)$ is isomorphic to an irreducible subquotient of $M_J((y_{x_2}x_2)\circ\lambda)$.
Then by Lemma~\ref{lem:subquotient} we have 
\begin{equation}
\label{eq:pr of th:4}
\ell(x_2)-\ell(y_{x_2})\leqq\ell(x_1)-\ell(y_{x_1}).
\end{equation}
On the other hand we have
\begin{equation}
\label{eq:pr of th:5}
\ell(x_2)+\ell(y_{x_2})=\ell(x_1)+\ell(y_{x_1})+1.
\end{equation}
by Lemma~\ref{lem:WJWI}.
Hence we have $2\ell(x_2)\leqq 2\ell(x_1)+1$.
Since $\ell(x_1)$ and $\ell(x_2)$ are integers, we obtain $\ell(x_2)\leqq \ell(x_1)$.
This is a contradiction.
The statement \eqref{eq:pr of th:1} is proved.

Let us show (i).
By Theorem~\ref{thm:main}(iv) we have $H^p(R(\D\O_{X_I}(\lambda)))=0$ for any $p\ne0$ if $\ell(x)\geqq m(x)$ for any $x\in\Gamma(\lambda)$.
Assume that $H^p(R(\D\O_{X_I}(\lambda)))=0$ for any $p>0$ and that $\ell(x)< m(x)$ for some $x\in\Gamma(\lambda)$.
Then we have $H^p(\shm(0)^\bullet)=0$ for any $p>0$ and $H^p(\shc(k)^\bullet)\ne0$ for some $p>0$ and some $k\geqq0$.
Let $r$ be the largest positive integer such that $H^r(\shc(k)^\bullet)\ne0$ for some $k\geqq0$.
Then we have an exact sequence
$$
H^r(\shm(k)^\bullet)\to H^r(\shm(k+1)^\bullet)\to 0\qquad(k\geqq0).
$$
Since  $H^r(\shm(0)^\bullet)=0$, we see by induction on $k$ that $H^r(\shm(k)^\bullet)=0$ for any $k\geqq0$.
Hence by \eqref{eq:pr of th:1} we have $H^r(\shc(k)^\bullet)=0$ for any $k\geqq0$.
This is a contradiction.
The statement (i) is proved.

Let us show (ii).
By (i) and the assumption we have $\ell(x)\geqq m(x)$ for any $x\in\Gamma(\lambda)$; in other words $H^p(\shc(k)^\bullet)=0$ for any $p>0$ and any $k\geqq0$.
By Theorem~\ref{thm:main}(v) $\Phi$ is an epimorphism if $\ell(x)> m(x)$ for any $x\in\Gamma(\lambda)\setminus\{e\}$.
Assume that $\Phi$ is an epimorphism.
Since $\Phi:H^0(\shc(0)^\bullet)\to H^0(\shm(0)^\bullet)$ is an epimorphism, we have $H^0(\shm(k)^\bullet)=0$ for any $k>0$ by the exact sequences
\begin{align*}
&H^0(\shc(0)^\bullet)\to H^0(\shm(0)^\bullet)\to H^0(\shm(1)^\bullet)\to 0,\\
&H^0(\shm(k)^\bullet)\to H^0(\shm(k+1)^\bullet)\to 0
\end{align*}
Hence by \eqref{eq:pr of th:1} we have $H^0(\shc(k)^\bullet)=0$ for any $k>0$.
It implies that  $\ell(x)> m(x)$ for any $x\in\Gamma(\lambda)\setminus\{e\}$.
The statement (ii) is proved.

Let us finally show (iii).
By (i) and the assumption we have $H^p(\shc(k)^\bullet)=0$ for any $p>0$ and any $k\geqq0$.
By Theorem~\ref{thm:main}(v) $\Phi$ is an isomorphism if $\Gamma(\lambda)=\{e\}$.
Hence it is sufficient to show that $H^{-p}(\shc(k)^\bullet)=0$ for any $k>0$ and any $p\geqq0$ if $\Phi$ is an isomorphism.
Let us show it by induction on $p$.
If $p=0$, then we have $H^0(\shc(k)^\bullet)=0$ for any $k>0$ by the proof of (ii).
Assume that the statement is proved up to $p$.
Consider the exact sequence
$$
H^{-(p+1)}(\shm(0)^\bullet)\to H^{-(p+1)}(\shm(1)^\bullet)\to H^{-p}(\shc(0)^\bullet)\to H^{-p}(\shm(0)^\bullet).
$$
We have $H^{-p}(\shc(0)^\bullet)=0$ for $p>0$, and $\Phi:H^{-p}(\shc(0)^\bullet)\to H^{-p}(\shm(0)^\bullet)$ is an isomorphism for $p=0$.
Moreover, we have $H^{-(p+1)}(\shm(0)^\bullet)=0$ by Lemma~\ref{lem:lower vanishing}.
Hence we have $H^{-(p+1)}(\shm(1)^\bullet)=0$.
Thus we obtain $H^{-(p+1)}(\shm(k)^\bullet)=0$ for any $k>0$ by the exact sequence
$$
H^{-(p+1)}(\shm(k)^\bullet)\to H^{-(p+1)}(\shm(k+1)^\bullet)\to H^{-p}(\shc(k)^\bullet)
$$
and the hypothesis of induction.
Hence we have $H^{-(p+1)}(\shc(k)^\bullet)=0$ for any $k>0$ by \eqref{eq:pr of th:1}.
The statement (iii) is proved.
\qed

By using Theorem~\ref{thm:extremal case} (i) and a case-by-case analysis we obtain the following.

\begin{theorem}
\label{thm:extremal vanishing}
Assume that $G$ is a simple group of classical type and that $\sharp(I)=\sharp(J)=\sharp(I_0)-1$.
Let $(\lambda,\mu)$ be the extremal case.
Then we have 
$H^p(R(\D\O_{X_I}(\lambda)))=0$ unless $p=0$.
\end{theorem}
Details of the proof is omitted.
We do not know an example of an extremal case $(\lambda,\mu)$ satisfying $H^p(R(\D\O_{X_I}(\lambda)))\ne0$ for some $p>0$.
\begin{remark}
Let $(\lambda,\mu)$ be an extremal case.
For $x\in\Gamma$ and $\alpha\in\Delta_J^+\setminus\Delta_I$ we have
$$
(x(\lambda+\rho))(\alpha^\vee)
=(\lambda+x\rho))(\alpha^\vee)
=(\mu-\gamma_{I,J}+x\rho)(\alpha^\vee)
=(x\rho-\gamma_{I,J})(\alpha^\vee),
$$
and hence we have $H^p(R(\D\O_{X_I}(\lambda)))=0$ for any $p>0$ if and only if
\begin{equation}
\label{eq:vanishing condition2}
\left\{\begin{array}{ll}
\mbox{for $x\in\Gamma$ satisfying $(x\rho-\gamma_{I,J})(\alpha^\vee)\ne0$ for any $\alpha\in\Delta_J^+\setminus\Delta_I$}\\
\mbox{we have $\sharp\{\alpha\in\Delta_J^+\setminus\Delta_I\,:\,
(x\rho-\gamma_{I,J})(\alpha^\vee)>0\}\leqq\ell(x)$ .}
\end{array}
\right.
\end{equation}
\end{remark}
In the next section we shall give conditions in order that $\Phi$ is an epimorphism and that $\Phi$ is an isomorphism under the assumption of Theorem~\ref{thm:extremal vanishing}.
In particular, $\Phi$ is not necessarily an epimorphism nor a monomorphism.
It seems to be an interesting problem to determine the kernel and the cokernel of $\Phi$.
\begin{remark}
In Tanisaki~\cite{Ta} it is shown in certain cases that $\Ker\Phi$ corresponds to the unique maximal proper submodule of $M_J(\mu)$ under the category equivalence given in Proposition~\ref{pr:G/H}.
\end{remark}

\section{The maximal parabolic case for classical simple groups}

In this section we apply our results to the case where $G$ is a classical simple group and $P_I$, $P_J$ are maximal parabolic subgroups, and obtain results for the Radon transform $R(\D\O_{X_I}(\lambda))$ with respect to the geometric correspondence  
$$
X_I\overset{f}{\longleftarrow} X_{I\cap J}\overset{g}{\longrightarrow}X_J
$$
for $\lambda\in(\gothh^*_\Z)^0_I$.
In this case we have
\begin{equation}
\mbox{$I=I_{0}\setminus\{p\}$ and $J=I_{0}\setminus\{q\}$ for some $p\neq q$,}
\end{equation}
and $(\gothh^*_\Z)^0_I=\{r\varpi_p:r\in\Z\}$, where $\varpi_k$ denotes the fundamental weight corresponding to $k\in I_0$.

We keep the standard notations of Bourbaki \cite{Bour}.
In particular, if $G$ is of rank $n$, then $I_0=\{1, 2, \dots, n\}$.

\subsection{The case (${\bf A_{n}}$)}
In this subsection we consider the case where $G=SL(V)$ for an $n+1$-dimensional complex vector space $V$.
By the symmetry of the Dynkin diagram we may (and shall) assume that $p>q$.
We have the identifications:
\begin{align*}
X_I&=\{\mbox{$p$-dimensional subspace of $V$}\},\\
X_J&=\{\mbox{$q$-dimensional subspace of $V$}\},\\
X_{I\cap J}&=
\{(U_1,U_2)\in X_I\times X_J:U_1\supset U_2\},
\end{align*}
and $f$, $g$ are natural projections.
The invertible $\O_{X_I}$-module $\O_{X_I}(\varpi_p)$ corresponds to the tautological line bundle whose fiber at $U\in X_I$ is $\bigwedge^pU$ (a subbundle of the product bundle $X_I\times\bigwedge^pV$), and we have $\O_{X_I}(r\varpi_p)=\O_{X_I}(\varpi_p)^{\otimes r}$.
Hence in the standard notation of algebraic geometry we have $\O_{X_I}(r\varpi_p)=\O_{X_I}(-r)$.

For $k\in I_0=\{1,\dots, n\}$
set 
$$
k_*=n+1-k,\quad
k_+=\max\{k,k_*\},\quad
k_-=\min\{k,k_*\}.
$$

We first give consequences of Theorem~\ref{thm:main}.
A weight $\lambda=\sum_{i=1}^{n+1} \lambda_i\eps_i$ ($\lambda_i\in\Z$, $\sum_{i=1}^{n+1}\lambda_i=0$) belongs to $(\gothh^*_\Z)_J$ if and only if 
$\lambda_1\geqq \cdots \geqq \lambda_q$ and $\lambda_{q+1}\geqq \cdots \geqq\lambda_{n+1}$.
The Weyl group $W$ is identified with the symmetric group $S_{n+1}$, and it acts on the weights by permutations of the components, i.\ e.\ $\sigma\lambda=\sum_{i=1}^{n+1}\lambda_i\eps_{\sigma(i)}$ for any $\sigma\in W$. 
Then we have $W_I=S_p\times S_{p_*}$ and $W_J=S_q\times S_{q_*}$.
We have
\begin{align*}
\varpi_p&=\frac1{n+1}\left[(n+1-p)(\eps_1+\cdots\eps_p)-p(\eps_{p+1}+\cdots\eps_{n+1})\right]\\
&=\eps_1+\cdots+\eps_p+\text{const}(\eps_1+\cdots+\eps_{n+1})\\
\rho&=\frac12\left[n\eps_1+(n-2)\eps_2+\cdots+(-n)\eps_{n+1}\right]\\
&=-\eps_2-\cdots-n\eps_{n+1}+\text{const}(\eps_1+\cdots+\eps_{n+1}),
\end{align*}
and therefore we get 
\begin{align*}
r\varpi_p+\rho&=r\eps_1+(-1+r)\eps_2+\cdots+(-(p-1)+r)\eps_p -p\eps_{p+1}
-\cdots -n\eps_{n+1}\\
&\qquad+\text{const}(\eps_1+\cdots+\eps_{n+1}).
\end{align*}
By the assumption $q<p$ the set $\Gamma(r\varpi_p)$ consists of $(\sigma,\tau)\in S_p\times S_{p_*}$ satisfying
$$
\left\{
\begin{array}{l}
\tau=e,\\
\sigma^{-1}(1)<\cdots<\sigma^{-1}(q),\\
\sigma^{-1}(q+1)<\cdots<\sigma^{-1}(p),\\
\{\sigma^{-1}(q+1),\dots,\sigma^{-1}(p)\}\cap\{p+r+1,\dots,n+r+1\}=\emptyset,
\end{array}
\right.
$$
and we have
\begin{align*}
\ell((\sigma,e))=&
\sharp\{(a,b)\,:\,1\leqq a\leqq q,\,q+1\leqq b\leqq p,\,
\sigma^{-1}(a)>\sigma^{-1}(b)\},\\
m((\sigma,e))=&
\sharp\{(b,c)\,:\,q+1\leqq b\leqq p,\,p+1\leqq c\leqq n+1,\,
\sigma^{-1}(b)<r+c\}.
\end{align*}
Hence by Theorem~\ref{thm:main} we obtain the following results.
\begin{proposition}
\label{pr:A:main}
\begin{itemize}
\item[\rm(i)]
Assume $q<p_-$.
Then we have $R(\D\O_{X_I}(-a\varpi_{p}))=0$ if $q+1\leqq a\leqq q_*-1$.
\item[\rm(ii)]
Assume $q\leqq p_-$.
Then we have 
\begin{align*}
R(\D\O_{X_I}(-q_*\varpi_{p}))=&\D\O_{X_J}(-p_*\varpi_{q}),\\
R(\D\O_{X_I}(-q\varpi_{p}))=&\D\O_{X_J}(-p\varpi_{q})[-(p-q)(p_*-q)].
\end{align*}
\item[\rm(iii)]
$H^k(R(\D\O_{X_I}(-a\varpi_p)))=0$ for any $k\ne0$ if $a\geqq q_*$.
\end{itemize}
\end{proposition}
Let us consider the extremal case.
By
$$
\gamma_{I,J}=p_*\sum_{i=q+1}^p\eps_i-(p-q)\sum_{i=p+1}^{n+1}\eps_i.
$$
and \eqref{eq:extrsys} the extremal case is given by $(-q_*\varpi_p,-p_*\varpi_q)$.
By Theorem~\ref{thm:extremal case} we obtain the following.

\begin{proposition}
\label{pr:A:extremal case}
We have $H^k((R(\D\O_{X_I}(-q_*\varpi_p)))=0$ for any $k\ne0$, and there exists a canonical nontrivial epimorphism
$$
\Phi:\D\O_{X_J}(-p_*\varpi_q)\to H^0(R(\D\O_{X_I}(-q_*\varpi_p))).
$$
Moreover, $\Phi$ is an isomorphism if and only if $p^*\geqq q$.
\end{proposition}
\begin{remark}
In the situation of Proposition~\ref{pr:A:extremal case} it is proved in \cite{Ta} that for $p^*\leqq q$ the kernel of $\Phi$ is the maximal proper $G$-stable submodule of $\D\O_{X_J}(-p_*\varpi_q)$.
\end{remark}

In the rest of this subsection we assume that $q<p_-$ and give application to topological problems.
By Proposition~\ref{pr:A:main} we have
\begin{equation}
\label{eq:SLDmod} R(\D\O_{X_I}(-a\varpi_{p})) \simeq
\left\{\begin{array}{ll}
\D\O_{X_J}(-p_*\varpi_{q}) & \mbox{for $a=q_*$,} \\ 
0 & \mbox{for $q+1\leq a\leq q_*-1$,} \\
\D\O_{X_J}(-p\varpi_{q})[-l_{pq}] & \mbox{for $a=q$,}
\end{array}\right.
\end{equation} 
where $l_{pq}=(p-q)(p_*-q)$.
Thus by Proposition~\ref{pr:adjunction} we have the following.

\begin{proposition}
\label{pr:A:adjunction} For any
$F\in\BDC(\C_{X_{J,\an}})$  and $q+1\leqq a\leqq q_*-1$ 
we have
\begin{eqnarray*}
\rsect(X_{I,\an};r(F)\tens\O_{X_I}(a\varpi_p)_\an) &=& 0,
\\
\RHom(r(F),\O_{X_I}(a\varpi_p)_\an) &=& 0,
\end{eqnarray*}
and for $(a,b,c,d)=(q_*,p_*,(p-q)p_*,pp_*-qq_*-q(p-q))$ 
or $(a,b,c,d)=(q,p,q(p-q),-q(p-q))$ we have
\begin{eqnarray*}
\rsect({X_{I,\an}};r(F)\tens\O_{X_I}(a\varpi_p)_\an)
&\simeq&
\rsect({X_{J,\an}};F\tens\O_{X_J}(b\varpi_q)_\an)[-c],
\\ 
\RHom(r(F),\O_{X_I}(a\varpi_p)_\an)
&\simeq&
\RHom(F,\O_{X_J}(b\varpi_q)_\an)[-d].
\end{eqnarray*}
\end{proposition}
Let us treat some particular cases. In the following
we set $N=qq_*$.
\vskip .3cm
\noindent 
(1) Let $y_\circ\in {X_J}$, and set $F=\C_{\{y_\circ\}}$.
Since $\opb{g}(y_{\circ})\to X_{I,y_\circ}$ is a closed
embedding, one has
\begin{equation}
\label{eq:rFpt} 
r(F)\simeq\C_{X_{I,y_\circ}{}_{,\an}},
\end{equation} 
where $X_{I,y_\circ}=f\opb{g}(\{y_\circ\})=\{x\in
{X_I}: y_\circ\subset x\}$ (identified with the Grassmannian
of $(p-q)$-subspaces of $V/y_\circ$). By
Proposition~\ref{pr:A:adjunction} and \eqref{eq:rFpt}
we obtain the following.
\begin{proposition} 
For any $q+1\leqq a\leqq q_*-1$ we have
$$
\rsect(X_{I,y_\circ}{}_{,\an};\O_{X_I}(a\varpi_p)_\an)
\simeq 0,
\qquad
\rsect_{X_{I,y_\circ}{}_{,\an}}({X_I};\O_{X_{I}}
(a\varpi_p)_\an)
\simeq 0,
$$ 
and for $(a,c,d)=(q_*,(p-q)p_*,pp_*-qq_*+p_*q)$ 
or $(a,c,d)=(q,q(p-q),p_*q)$ we have
$$
H^{c}(X_{I,y_\circ}{}_{,\an};\O_{X_I}(a\varpi_p)_\an)
\simeq \C\{z\},
\qquad
H^{d}_{X_{I,y_\circ}{}_{,\an}}({X_I};\O_{X_{I}}(a\varpi_p)_\an)
\simeq
\shb^{\infty}_{0|\C^N}
$$
where $\C\{z\}$ (resp.\ $\shb^{\infty}_{0|\C^N}$) is the
ring of convergent power series in $z=(z_1,\dots,z_N)\in\C^N$
(resp.\ the ring of hyperfunctions in $\C^N$ along $\{0\}$
of infinite order), and all other cohomology groups vanish.
\end{proposition}

\noindent
Namely, one identifies $\rsect(X_{J,\an};\C_{y_\circ}\tens\O_{X_J}
(b\varpi_q)_\an)\simeq\rsect(\{0\};\O_{\C^N_\an})=\C\{z\}$ and
$\RHom(\C_{y_\circ};\O_{X_J}(b\varpi_q)_\an)\simeq
\rsect_{\{0\}}(\C^N_\an;\O_{\C^N_\an})=\shb^{\infty}_{0|\C^N}[-N]$.
\vskip .3cm
\noindent 
(2) Let $z_{\circ}$ be a $q_*$-subspace of $V$,
$E_{z_{\circ}}=\{y\in{X_J}:y\cap z_{\circ}=0\}
\simeq\C^N$ and set $F=\C_{E_{z_{\circ},\an}}$. 
One has
\begin{equation}
\label{eq:rF:hyperplane} 
r(F)\simeq\C_{\widehat{E}_{z_{\circ}}{}_{,\an}}[-2q(p-q)],
\end{equation} 
where $\widehat{E_{z_{\circ}}}=f\opb{g}(E_{z_{\circ}})=
\{x\in {X_I}: \dim(x\cap{z_{\circ}})=p-q\}$ (i.e.\ the 
$p$-dimensional subspaces of $V$ in generic position 
w.r.t.\ ${z_{\circ}}$). Namely, the map
$\tilde{f}=(\left.f\right|_{\opb{g}(E_{z_{\circ}})})_\an:
(\opb{g}(E_{z_{\circ}}))_\an
\to \widehat{E}_{z_{\circ}}{}_{,\an}$ 
is a complex vector bundle of rank $q(p-q)$
(the fiber over $x\in \widehat{E}_{z_{\circ}}$ is
$S_{E_{z_{\circ}},x}=\{y\in E_{z_{\circ}}:y\subset
x\}\simeq\C^{q(p-q)}$); hence there is a  morphism of
functors $\roim{\tilde{f}}\opb{\tilde{f}}[2q(p-q)]\to
\id_{\BDC(\C_{\widehat{E_{z_{\circ}}}{}_{,\an}})}$
defining a natural morphism $r(F)=\reim{f_\an} 
\C_{(\opb{g}(E_{z_{\circ}}))_\an}
\to\C_{\widehat{E}_{z_{\circ}}{}_{,\an}}[-2q(p-q)]$,
which is an isomorphism since, by \eqref{eq:rfiber},
one has $r(F)_x\simeq \C[-2q(p-q)]$ (for 
$x\in\widehat{E}_{z_{\circ}}{}_{,\an}$) and $=0$
(otherwise).

By Proposition~\ref{pr:A:adjunction} and
\eqref{eq:rF:hyperplane} we obtain the following.
\begin{proposition} 
For any $q+1\leqq a\leqq q_*-1$ we have
$$
\rsect_{\rm c}(\widehat{E}_{z_{\circ}}{}_{,\an};
\O_{X_I}(a\varpi_p)_\an)\simeq 0,
\qquad
\rsect(\widehat{E}_{z_{\circ}}{}_{,\an};
\O_{X_I}(a\varpi_p)_\an)\simeq 0,
$$ 
and for $(a,c,d)=(q_*,p(p_*-q)+q^2,p_*(p-q))$ 
or $(a,c,d)=(q,p_*q,q(p-q))$ we have
\begin{eqnarray*}
H^{c}_{\rm c}(\widehat{E}_{z_{\circ}}{}_{,\an};
\O_{X_I}(a\varpi_p)_\an)
&\simeq& H^N_{\rm c}(E_{z_{\circ}}{}_{,\an};\O_{E_{z_{\circ}}
{}_{,\an}}),
\\
H^{d}(\widehat{E}_{z_{\circ}}{}_{,\an};
\O_{X_I}(a\varpi_p)_\an) &\simeq&
\sect(E_{z_{\circ}}{}_{,\an};\O_{E_{z_{\circ}}{}_{,\an}})
\end{eqnarray*}
where $H^N_{\rm c}(E_{z_{\circ}}{}_{,\an};
\O_{E_{z_{\circ}}{}_{,\an}})
\simeq \sect(E_{z_{\circ}}{}_{,\an};
\Omega_{E_{z_{\circ}}{}_{,\an}})'$
(resp.\ $\sect(E_{z_{\circ}}{}_{,\an};
\O_{E_{z_{\circ}}{}_{,\an}})$)
are Martineau's analytic functionals (resp.\ the entire
functions) in $E_{z_{\circ}}{}_{,\an}\simeq\C^N$, and all
other cohomology groups vanish.
\end{proposition}
Namely, one has
$\rsect(X_{J,\an};\C_{E_{z_{\circ}}{}_{,\an}}\tens\O_{X_J}
(b\varpi_q)_\an)\simeq
H^N_{\rm c}(E_{z_{\circ}}{}_{,\an};
\O_{E_{z_{\circ}}{}_{,\an}})[-N]$ 
and $\RHom(\C_{E_{z_{\circ}}{}_{,\an}};
\O_{X_J}(b\varpi_q)_\an) \simeq
\sect(E_{z_{\circ}}{}_{,\an};
\O_{E_{z_{\circ}}{}_{,\an}})$.

\subsection{The case (${\bf B_{n}}$)}
In this subsection we consider the case where $G$ is (the universal covering group of) $SO(V)$ for an $2n+1$-dimensional complex vector space $V$ equipped with a non-degenerate symmetric bilinear form $(\,,\,):V\times V\to\C$.
Then we have the identifications:
\begin{align*}
X_I&=
\{\mbox{$p$-dimensional subspace $U$ of $V$ such that $(U,U)=0$}\},\\
X_J&=
\{\mbox{$q$-dimensional subspace $U$ of $V$ such that $(U,U)=0$}\},\\
X_{I\cap J}&=
\left\{
\begin{array}{ll}
\{(U_1,U_2)\in X_I\times X_J:U_1\subset U_2\}\qquad&(p<q)\\
\{(U_1,U_2)\in X_I\times X_J:U_1\supset U_2\}\qquad&(p>q),
\end{array}
\right.
\end{align*}
and $f$, $g$ are natural projections.
The invertible $\O_{X_I}$-module $\O_{X_I}(\varpi_p)$ corresponds to the tautological line bundle whose fiber at $U\in X_I$ is $\bigwedge^pU$.

By Theorem~\ref{thm:main} we have the following.
\begin{proposition}
\label{pr:B:main}
\begin{itemize}
\item[\rm(i)]
We have $R(\D\O_{X_I}(-a\varpi_p))=0$ 
in the following cases:
$$
\left\{
\begin{array}{ll}
2n-p-q<a<q \qquad&\mbox{if $p<q\leqq n$,}\\
\min(n-p,q)<a<n-q  \qquad&\mbox{if $q<p<n$,}\\
n-p+q<a<\max(n,2n-p-q)  \qquad&\mbox{if $q<p<n$,}\\
2q<a<2(n-q)\qquad&\mbox{if $p=n$,}\\
\end{array}
\right.
$$
\item[\rm(ii)]
We have $R(\D\O_{X_I}(-a\varpi_p))=\D\O_{X_J}(-b\varpi_q)[-c]$ in the following cases:
\begin{align*}
&(a, b, c)\\
=&\left\{
 \begin{array}{ll}
  (q,p,0)\quad&
  (p<q<n,\, 2n-2p-q\leqq0),\\
  (2n-p-q, 2(n-q), c_1)\quad&
  (p<q\leqq n,\, 2n-2p-q\leqq0),\\
  (n, 2p,0)\quad&
  (q=n, \, n-2p\leqq0),\\
  (2n-p-q,2n-p-q,0)\quad&
  (q<p<n, \, 2n-2p-q\geqq0),\\
  (q,p,c_2)\quad&
  (q<p<n, \, 2n-2p-q\geqq0),
 \end{array}
\right.
\end{align*}
where
\[
c_1=\frac{(q-p)(3p+3q-4n-1)}{2},\qquad
c_2=\frac{(p-q)(4n+1-3p-3q)}{2}.
\]
\end{itemize}
\end{proposition}
By Theorem~\ref{thm:extremal case} we have the following.
\begin{proposition}
\label{pr:B:extremal case}
Let
$$
(r,s)=\left\{
\begin{array}{ll}
(q,p)\qquad&\mbox{if $1\leqq p<q\leqq n-1$,}\\
(2n-p-q,2n-p-q)\qquad&\mbox{if $1\leqq q< p\leqq n-1$,}\\
(2(n-q),n-q)\qquad&\mbox{if $p=n$, $1\leqq q\leqq n-1$,}\\
(n,2p)\qquad&\mbox{if $1\leqq p\leqq n-1$, $q=n$.}
\end{array}
\right.
$$
Then we have $H^k(R(\D\O_{X_I}(-r\varpi_p)))=0$ for any $k\ne0$, and there exists a canonical nontrivial morphism
$$
\Phi:\D\O_{X_J}(-s\varpi_q)\to H^0(R(\D\O_{X_I}(-r\varpi_p))).
$$
Moreover, $\Phi$ is an epimorphism if and only if we have either 
\begin{itemize}
\item[\rm(a)]
$p<q\leqq n$, 
\item[\rm(b)]
$q<p<n$ and $2n-2p-q\geqq0$,
\end{itemize}
and an isomorphism if and only if we have either
\begin{itemize}
\item[\rm(a)]
$p<q\leqq n$ and $2n-2p-q\leqq0$, 
\item[\rm(b)]
$q<p<n$ and $2n-2p-q\geqq0$.
\end{itemize}
\end{proposition}

\subsection{The case (${\bf C_{n}}$)}
In this subsection we consider the case where $G=Sp(V)$ for an $2n$-dimensional complex vector space $V$ equipped with a non-degenerate anti-symmetric bilinear form $(\,,\,):V\times V\to\C$.
Then we have the identifications:
\begin{align*}
X_I&=
\{\mbox{$p$-dimensional subspace $U$ of $V$ such that $(U,U)=0$}\},\\
X_J&=
\{\mbox{$q$-dimensional subspace $U$ of $V$ such that $(U,U)=0$}\},\\
X_{I\cap J}&=
\left\{
\begin{array}{ll}
\{(U_1,U_2)\in X_I\times X_J:U_1\subset U_2\}\qquad&(p<q)\\
\{(U_1,U_2)\in X_I\times X_J:U_1\supset U_2\}\qquad&(p>q),
\end{array}
\right.
\end{align*}
and $f$, $g$ are natural projections.
The invertible $\O_{X_I}$-module $\O_{X_I}(\varpi_p)$ corresponds to the tautological line bundle whose fiber at $U\in X_I$ is $\bigwedge^pU$.

By Theorem~\ref{thm:main} we have the following.
\begin{proposition}
\label{pr:C:main}
\begin{itemize}
\item[\rm(i)]
We have $R(\D\O_{X_I}(-a\varpi_p))=0$ 
in the following cases:
$$
\left\{
\begin{array}{ll}
2n-p-q+1<a<q \qquad&\mbox{if \, $p<q$,}\\
q<a<2n-p-q+1 \qquad&\mbox{if \, $q< p$.}
\end{array}
\right.
$$
\item[\rm(ii)]
We have $R(\D\O_{X_I}(-a\varpi_p))=\D\O_{X_J}(-b\varpi_q)[-c]$ in the following cases:
\begin{align*}
&(a, b, c)\\
=&\left\{
 \begin{array}{ll}
  (q,p,0)\quad&
  (p<q\leqq n,\, 2n-2p-q+1\leqq0),\\
  (2n-p-q+1, 2n-p-q+1, c_1)\quad&
  (p<q\leqq n,\, 2n-2p-q+1\leqq0),\\
  (2n-p-q+1,2n-p-q+1,0)\quad&
  (q<p\leqq n, \, 2n-2p-q+1\geqq0),\\
  (q,p,c_2)\quad&
  (q<p\leqq n, \, 2n-2p-q+1\geqq0),
 \end{array}
\right.
\end{align*}
where
\[
c_1=\frac{(q-p)(3p+3q-4n-1)}{2},\qquad
c_2=\frac{(p-q)(4n+1-3p-3q)}{2}.
\]
\end{itemize}
\end{proposition}

By Theorem~\ref{thm:extremal case} we have the following.
\begin{proposition}
\label{pr:C:extremal case}
Let
$$
(r,s)=\left\{
\begin{array}{ll}
(q,p)\qquad&\mbox{if $1\leqq p<q\leqq n$,}\\
(2n-p-q+1,2n-p-q+1)\qquad&\mbox{if $1\leqq q< p\leqq n$.}
\end{array}
\right.
$$
Then we have $H^k(R(\D\O_{X_I}(-r\varpi_p)))=0$ for any $k\ne0$, and there exists a canonical nontrivial morphism
$$
\Phi:\D\O_{X_J}(-s\varpi_q)\to H^0(R(\D\O_{X_I}(-r\varpi_p))).
$$
Moreover, $\Phi$ is an epimorphism if and only if we have either 
\begin{itemize}
\item[\rm(a)]
$p<q<n$ and $n-p-q\geqq0$,
\item[\rm(b)]
$p<q\leqq n$ and $2n-2p-q+1\leqq 0$,
\item[\rm(c)]
$q<p\leqq n$,
\end{itemize}
and an isomorphism if and only if we have either
\begin{itemize}
\item[\rm(a)]
$p<q\leqq n$ and $2n-2p-q+1\leqq 0$,
\item[\rm(b)]
$q<p\leqq n$ and $2n-2p-q+1\geqq 0$.
\end{itemize}
\end{proposition}
\begin{remark}
In the situation of Proposition~\ref{pr:C:extremal case} it is proved in \cite{Ta} that $\Ker\Phi$ is the maximal proper $G$-stable submodule of $\D\O_{X_J}(-s\varpi_q)$ if $q=n$ and $2p\leqq n-1$.
\end{remark}

\subsection{The case (${\bf D_{n}}$)}
In this subsection we consider the case where $G$ is (the universal covering group of) $SO(V)$ for an $2n$-dimensional complex vector space $V$ equipped with a non-degenerate symmetric bilinear form $(\,,\,):V\times V\to\C$.

For $1\leqq k\leqq n$ set
$$
X(k)=\{\mbox{$k$-dimensional subspace $U$ of $V$ such that $(U,U)=0$}\}.
$$
Then $X(k)$ is connected for $1\leqq k\leqq n-1$, and $X(n)$ has two connected components, say $X_1(n)$ and $X_2(n)$.
Then we have the identification:
\begin{align*}
&X(k)=X_{I_0\setminus\{k\}}\qquad(1\leqq k\leqq n-2),\\
&X(n-1)=X_{I_0\setminus\{n-1, n\}},\\
&X_1(n)=X_{I_0\setminus\{n\}},\\
&X_2(n)=X_{I_0\setminus\{n-1\}}.
\end{align*}
If $\{p, q\}\ne\{n-1, n\}$, then
$$
X_{I\cap J}=
\left\{
\begin{array}{ll}
\{(U_1,U_2)\in X_I\times X_J:U_1\subset U_2\}\qquad&(p<q)\\
\{(U_1,U_2)\in X_I\times X_J:U_1\supset U_2\}\qquad&(p>q),
\end{array}
\right.
$$
and if $p=n-1$ and $q=n$, then $f$ (resp.\ $g$) assigns $U\in X_{I\cap J}=X(n-1)$ to the unique $U'\in X_I=X_2(n)$ (resp.\ $U'\in X_J=X_1(n)$) such that $U\subset U'$.
The invertible $\O_{X_I}$-module $\O_{X_I}(\varpi_p)$ corresponds to the tautological line bundle whose fiber at $U\in X_I$ is $\bigwedge^kU$ where $k=p$ for $1\leqq k\leqq n-2$ and $k=n$ for $p\in\{n-1, n\}$.

By Theorem~\ref{thm:main} we have the following.
\begin{proposition}
\label{pr:D:main}
\begin{itemize}
\item[\rm(i)]
We have $R(\D\O_{X_I}(-a\varpi_p))=0$ 
in the following cases:
$$
\left\{
\begin{array}{ll}
2n-p-q-1<a<q\quad&\mbox{if $p<q\leqq n-2$,}\\
q<a<2n-p-q-1\quad&\mbox{if $q< p\leqq n-2$}\\
2q<a<2(n-q-1)\quad&\mbox{if $p\in\{n-1,n\}$, $1\leqq q\leqq n-2$,}\\
n-p-1<a<n\quad&\mbox{if $1\leqq p\leqq n-2$, $q\in\{n-1,n\}$,}\\
a=n-1\quad&\mbox{if $\{p,q\}=\{n-1,n\}$ and $n$ is even.}
\end{array}
\right.
$$
\item[\rm(ii)]
We have $R(\D\O_{X_I}(-a\varpi_p))=\D\O_{X_J}(-b\varpi_q)[-c]$ in the following cases:
{\small
\begin{align*}
&(a, b, c)\\
=&\left\{
 \begin{array}{ll}
  (q,p,0)\quad&
  (p<q\leqq n-2,\, 2n-2p-q-1\leqq0),\\
  (2n-p-q-1, 2n-p-q-1, c_1)\quad&
  (p<q\leqq n-2,\, 2n-2p-q-1\leqq0),\\
  (n,2p,0)\quad&
  (p\leqq n-2,\, q\in\{n-1,n\}, \, n-2p-1\leqq0),\\
  (n-p-1, 2(n-p-1), c_2)\quad&
  (p\leqq n-2,\, q\in\{n-1,n\}, \, n-2p-1\leqq0),\\
  (2n-p-q-1,2n-p-q-1,0)\quad&
  (q<p\leqq n-2, \, 2n-2p-q-1\geqq0),\\
  (q,p,c_3)\quad&
  (q<p\leqq n-2, \, 2n-2p-q-1\geqq0),\\
  (n,n-2,0)\quad&
  (\{p,q\}=\{n-1,n\}, \, n:\mbox{odd}),\\
  (n-1,n-1,0)\quad&
  (\{p,q\}=\{n-1,n\}, \, n:\mbox{odd}),\\
  (n-2,n,0)\quad&
  (\{p,q\}=\{n-1,n\}, \, n:\mbox{odd}),
  \end{array}
\right.
\end{align*}
}
where 
\begin{align*}
&c_1=\frac{(q-p)(3p+3q-4n+1)}{2},\qquad
c_2=\frac{(n-p)(3p-n+1)}{2},\\
&c_3=\frac{(p-q)(4n-3p-3q-1)}{2}.
\end{align*}
\end{itemize}
\end{proposition}

By Theorem~\ref{thm:extremal case} we have the following.
\begin{proposition}
\label{pr:D:extremal case}
Let
\begin{align*}
&(r,s)\\
=&\left\{
\begin{array}{ll}
(q,p)\quad&\mbox{if $1\leqq p<q\leqq n-2$,}\\
(2n-p-q-1,2n-p-q-1)\quad&\mbox{if $1\leqq q< p\leqq n-2$,}\\
(2(n-q-1),n-q-1)\quad&\mbox{if $p\in\{n-1,n\}$, $1\leqq q\leqq n-2$,}\\
(n,2p)\quad&\mbox{if $1\leqq p\leqq n-2$, $q\in\{n-1,n\}$,}\\
(n,n-2)\quad&\mbox{if $\{p,q\}=\{n-1,n\}$.}
\end{array}
\right.
\end{align*}
Then we have $H^k(R(\D\O_{X_I}(-r\varpi_p)))=0$ for any $k\ne0$, and there exists a canonical nontrivial epimorphism
$$
\Phi:\D\O_{X_J}(-s\varpi_q)\to H^0(R(\D\O_{X_I}(-r\varpi_p))).
$$
Moreover, $\Phi$ is an isomorphism if and only if we have either
\begin{itemize}
\item[\rm(a)]
$p<q<n-1$ and $2n-2p-q-1\leqq0$,
\item[\rm(b)]
$q<p<n-1$ and $2n-2p-q-1\geqq0$,
\item[\rm(c)]
$p<n-1$, $q\in\{n-1, n\}$ and $n-2p-1\leqq0$,
\item[\rm(c)]
$\{p, q\}=\{n-1, n\}$ and $n$ is odd.
\end{itemize}
\end{proposition}

\begin{remark}
In the situation of Proposition~\ref{pr:D:extremal case} it is proved in \cite{Ta} that $\Ker\Phi$ is the maximal proper $G$-stable submodule of $\D\O_{X_J}(-s\varpi_q)$ if 
$q\in\{n-1, n\}, 2p\leqq n-2$ and if $q=1, p\in\{n-1, n\}$.
\end{remark}

\vfill\break

\noindent
\underbar{\hphantom{XXXXXXXXXXXXXXXX}}

{\small

\vskip .3cm
\noindent
Corrado Marastoni: Universit\`a di Padova -- Dipartimento di Matematica Pura ed
Applicata -- Via Belzoni, 7 -- I-35131 Padova (Italy) -- {\tt
maraston@math.unipd.it}

\vskip .1cm
\noindent
Toshiyuki Tanisaki: Hiroshima University -- Department of Mathematics --
Higashi-Hiro\-shi\-ma 739-8526 (Japan) -- {\tt
tanisaki@math.sci.hiroshima-u.ac.jp}

}

\end{document}